\documentclass{siamart251216}

\usepackage{microtype}
\usepackage{amsfonts,amsmath,amssymb}
\usepackage{bm}
\usepackage{mathrsfs}
\usepackage{enumitem}
\usepackage{cases}
\usepackage{graphicx}
\usepackage{stmaryrd}
\usepackage{xparse}
\usepackage{epstopdf}
\usepackage{cite}
\usepackage{algorithmic}
\usepackage{subcaption}

\usepackage{multirow}
\graphicspath{{fig/}{fig/2d/}{fig/3d}}
\usepackage{color}
\ifpdf
  \DeclareGraphicsExtensions{.eps,.pdf,.png,.jpg}
\else
  \DeclareGraphicsExtensions{.eps}
\fi
\allowdisplaybreaks[4]

\newsiamremark{remark}{Remark}
\newsiamremark{numscheme}{Numerical Scheme}
\newsiamremark{assumption}{Assumption}
\newsiamremark{hypothesis}{Hypothesis}
\crefname{hypothesis}{Hypothesis}{Hypotheses}
\newsiamthm{claim}{Claim}

\def\Omega{\varOmega}
\def\Delta{\varDelta}
\def\R{{\mathbb R}}

\def\lj{{\bm [ }\!{\bm [}}
\def\rj{{\bm ] }\!{\bm ]}}

\def\eps{\epsilon}

\def\T{{\mathcal{T}}}

\headers{RFEM for the Monge--Amp\`ere equation
}{Tianyang Chu, Hailong Guo and Zhimin Zhang}

\title{A recovery based linear finite element method for the Monge--Amp\`{e}re equation}

\author{
Tianyang Chu\thanks{School of Mathematical Sciences, Chengdu University of Technology, Chengdu, Sichuan, 610059, China (\email{tchu@cdut.edu.cn}).}
\and Hailong Guo\thanks{Chongqing Research Institute of Big Data, Peking University, Chongqing, 401332, China (\email{hailong.guo@cqbdri.pku.edu.cn}).}
\and Zhimin Zhang\thanks{Department of Mathematics, Wayne State University, Detroit, MI 48202, USA (\email{ag7761@wayne.edu}).}
}

\begin{document}

\maketitle

\begin{abstract}
In this paper, we study finite element approximations of a fourth-order perturbed semilinear partial differential equation that arises naturally in the computation of viscosity solutions to the fully nonlinear Monge--Amp\`{e}re equation. We develop a novel recovery-based $C^0$ linear finite element method for discretizing this model problem. The key analytical tool is a Miranda--Talenti-type estimate for the recovery operator, which guides the formulation of the discrete scheme and guarantees its well-posedness. We rigorously establish convergence in the $H^2$ norm by employing a fixed-point argument. Furthermore, optimal $H^1$-error estimates are derived using a duality argument. Numerical experiments confirm the efficiency of the proposed method and validate the theoretical convergence results.
\end{abstract}



\section{Introduction}\label{sec: int}

In this paper, we consider the fully nonlinear Monge--Amp\`ere equation with Dirichlet boundary conditions:
\begin{subequations}\label{equ: MA equation}
    \begin{align}
        \det(D^2u) &= f \quad \text{in } \Omega, \\
        u &= g \quad \text{on } \partial\Omega,
    \end{align}
\end{subequations}
where $\Omega \subset \mathbb{R}^2$ is a convex polygonal domain with Lipschitz continuous boundary $\Gamma = \partial\Omega$, and $D^2u$ denotes the Hessian matrix of the scalar function $u$. The nonlinearity of the equation is captured by the determinant of the Hessian. Such equations arise in a variety of scientific and engineering contexts, such as in problems from differential geometry~\cite{guti2016ma, figalli2017mae}, optical reflector design~\cite{wang1996inversedesign}, and the theory of optimal transport~\cite{villani2003topicinot, benamou2000cfmot}. For the Monge--Amp\`ere equation, alternative notions of solutions, such as Aleksandrov and viscosity solutions~\cite{caffarelli1995fullynonlinear, figalli2017mae, guti2016ma, gilbarg2001secondorderpde}, have been developed over the past century, as classical solutions may not exist even when the data $f$ and $g$ are smooth.

Due to its nonlinear nature and lack of classical variational structure, the development of numerical techniques for the Monge--Amp\`ere equation is challenging.  The first numerical methods can be traced back to the late 1980s when Oliker and Prussner \cite{oliker1988geometric} constructed the finite difference methods to compute the Aleksandrov solutions.  In 2008, Oberman proposed a wide-stencil framework for solving the Monge--Amp\`ere equation \cite{oberman2008wide}. Thus far, a variety of numerical methods have been proposed for the approximation of Monge--Amp\`ere equations, including filtered schemes~\cite{froese2011confdm, froese2013filterscheme}, $L^2$-projection methods and $C^0$ interior penalty methods~\cite{bohmer2008l2projection, brenner2011c0ip2d, brenner2012c0ip3d, neilan2013quad, adetola2022recovery}, least-squares approaches~\cite{brenner2021convex, dean2004leastsquare}, vanishing moment methods~\cite{feng2009mixed, feng2009modified, neilan2010nonconforming, chenfengzhang2021, feng2009vanishing, feng2011conforming}, and two-scale finite element methods~\cite{nochetto2019twoscalepoint, nochetto2019twoscale}. For a comprehensive overview, we refer the reader to the review articles~\cite{neilan2017review, neilan2020review, feng2013siamreview}.

We concentrate on the finite element approximation of the following singularly perturbed fourth-order semilinear partial differential equation:
\begin{subequations}\label{equ: model equations}
    \begin{align}
        -\epsilon \varDelta^2 u^{\epsilon} + \det(D^2u^{\epsilon}) = f& \quad \text{in $\varOmega$}, \\
        u^{\epsilon} =g&\quad \text{on $\partial\varOmega$}, \\
        \varDelta u^{\epsilon} =\eps & \quad \text{on $\partial\varOmega$}.
    \end{align}
\end{subequations} 
The solution of the model problem~\eqref{equ: model equations} has been shown to converge uniformly (in the viscosity sense) to the viscosity solution of the Monge--Amp\`ere equation~\eqref{equ: MA equation} as $\epsilon \to 0^+$, in the radially symmetric setting~\cite{feng2014convergence, feng2009vanishing}. Thus, equation~\eqref{equ: model equations} offers a regularization framework for computing viscosity solutions of the Monge--Amp\`ere equation. Following the assumptions in~\cite{feng2009mixed, chenfengzhang2021}, we further suppose that the following a priori estimates hold:
\begin{equation}\label{equ: priori bounds for u^eps}
	\begin{aligned}
		&\|u^{\eps}\|_{H^{j}(\Omega)} = \mathcal{O}(\eps^{\frac{1-j}{2}}) \ (j = 2,\; 3), && \|u^{\eps}\|_{W^{j, \infty}(\Omega)} = \mathcal{O}(\eps^{1-j}) \ (j = 1,\; 2), \\
		&\|\Phi^{\eps}\|_{L^2(\Omega)} = \mathcal{O}(\eps^{-\frac{1}{2}}), && \|\Phi^{\eps}\|_{L^{\infty}(\Omega)} = \mathcal{O}(\eps^{-1}), 
	\end{aligned}
\end{equation}
where $\Phi^{\eps}$ denotes the cofactor matrix of the Hessian $D^2u^{\eps}$. Additionally, we assume that problem \eqref{equ: model equations} admits a unique convex solution $u^{\eps} \in H^3(\Omega) \cap W^{2,\infty}(\Omega)$, and satisfies \eqref{equ: priori bounds for u^eps}.

For the semilinear equation~\eqref{equ: model equations}, one of the primary numerical challenges lies in discretizing the biharmonic operator. Classical finite element methods for the biharmonic operator, such as conforming finite element methods~\cite{feng2011conforming}, mixed finite element methods~\cite{feng2009mixed}, nonconforming finite element methods~\cite{neilan2010nonconforming}, and $C^0$ interior penalty methods~\cite{neilan2020review, chu2024c0ip}, have been applied to solve~\eqref{equ: model equations}. Most of these studies use at least quadratic elements. In recent years, there has been increasing interest in designing linear finite element methods for solving fourth-order PDEs; see, e.g.,~\cite{Lami2014, larsson2012continuous, ChenGuoZhangZou2017, GuoZhangZou2018, XuGuoZou2019, cai2025stablized}. However, these approaches are not directly applicable to~\eqref{equ: model equations} due to difficulties in establishing stability conditions and the lack of consistency in the finite element functions. To the best of our knowledge,~\cite{chenfengzhang2021} is the first work devoted to the development of linear element methods for solving the model equation~\eqref{equ: model equations}. However, certain assumptions are required in~\cite{chenfengzhang2021} to ensure stability.

The primary goal of this paper is to propose a new linear element-based method with rigorous mathematical theory. The main motivation comes from our recent work on the Hamilton--Jacobi--Bellman equation~\cite{chu2025recovery}, where a novel gradient recovery operator is proposed using the Morley finite element function as a bridge. Unlike existing gradient recovery operators~\cite{guo2025recovery}, the proposed operator is a true gradient of a Morley finite element function, which allows us to establish the consistency of the gradient recovery operator for finite element functions. In addition, a discrete Miranda--Talenti (MT) estimate for the gradient recovery operator is established in~\cite{chu2025recovery}. Unlike the analysis method in~\cite{chenfengzhang2021}, our analysis uses the discrete MT estimate to design the finite element scheme and analyze the stability of the method. Based on the Banach fixed-point theorem, we prove the existence and uniqueness of the finite element solution, as well as a discrete $H^2$ error estimate. What we want to emphasize is that the lower-order norm error estimates based on the MT estimate are missing from the literature. To obtain the optimal error estimate in the $H^1$ norm, we employ a duality argument. Compared to the recovery-based finite element method in~\cite{chenfengzhang2021}, our method achieves $\mathcal{O}(h^2)$ ($\mathcal{O}(h)$) convergence in the $L^2$ norm ($H^1$ norm), while the method in~\cite{chenfengzhang2021} only yields $\mathcal{O}(h)$ ($\mathcal{O}(h^{0.5})$) convergence in the $L^2$ norm ($H^1$ norm).

The remainder of the paper is organized as follows. In Section~\ref{sec: pre},  we introduce the necessary notation for finite element functions and review the gradient recovery operator developed in~\cite{chu2025recovery}. Section~\ref{sec:lin} is devoted to the design and analysis of a finite element scheme for the linearized model, with particular emphasis on the application of a discrete MT estimate. In Section~\ref{sec:nonlinear}, we present a recovery-based finite element method for the semilinear equation and establish error estimates in the discrete $H^2$ and $H^1$ norms. Section~\ref{sec:num} contains several numerical experiments that illustrate the performance of the proposed method. Concluding remarks are provided in Section~\ref{sec:con}. Additional properties of the gradient recovery operator are discussed in Appendix~\ref{sec: appendix}. 

\section{Preliminaries}\label{sec: pre}

\subsection{Notation}\label{ssec:not}

Let $|D|$ denote the Lebesgue measure of a measurable set $D\subset \mathbb{R}^2$. For $s \geq 0$ and $p \in [1,\infty]$, we denote by $W^{s,p}(D)$ the standard Sobolev space on the domain $D$. In particular, we write $L^p(D) := W^{0,p}(D)$ and $H^s(D) := W^{s,2}(D)$. The corresponding norm and semi-norm on $W^{s,p}(D)$ are denoted by $\|\cdot\|_{W^{s,p}(D)}$ and $|\cdot|_{W^{s,p}(D)}$, respectively. Let $H_0^1(D)$ denote the subspace of $H^1(D)$ consisting of functions that vanish on $\partial D$ in the trace sense. The inner product in $L^2(D)$ is denoted by $(\cdot,\cdot)_{L^2(D)}$, and let 
\begin{equation*}
	V:= H^2(\Omega), \quad V_0:= H^2(\Omega)\cap H_0^1(\Omega), \quad V_g:= \{ v \in V: v|_{\partial\Omega}=g \}.
\end{equation*}

Let $\mathcal{T}_h$ be a shape-regular and conforming simplicial triangulation of the domain $\Omega$. For each element $T \in \mathcal{T}_h$, denote its diameter by $h_T := \mathrm{diam}(T)$, and let $h := \max_{T \in \mathcal{T}_h} h_T$ represent the mesh size. Let $\mathcal{V}_h$ and $\mathcal{E}_h$ denote the sets of vertices and edges in $\mathcal{T}_h$, respectively. For each edge $e \in \mathcal{E}_h$, let $h_e$ denote its length. The edge set $\mathcal{E}_h$ is partitioned into interior and boundary edges, denoted by $\mathcal{E}_h^I$ and $\mathcal{E}_h^B$, respectively. Likewise, the vertex set $\mathcal{V}_h$ is divided into interior and boundary vertices, denoted by $\mathcal{V}_h^I$ and $\mathcal{V}_h^B$. For each element $T \in \mathcal{T}_h$, we let $\mathcal{V}_T$ and $\mathcal{M}_T$ denote the sets of vertices and midpoints of the edges of $T$, respectively. We define $\mathcal{T}_p$ as the collection of elements in $\mathcal{T}_h$ that share a common vertex $p$, and $\mathcal{T}_e$ as the set of elements sharing a common edge $e \in \mathcal{E}_h$. Finally, let $\mathcal{M}_h^I$ and $\mathcal{M}_h^B$ denote the sets of interior and boundary edge midpoints, respectively, and define $\mathcal{M}_h := \mathcal{M}_h^I \cup \mathcal{M}_h^B$.

For any edge $e \in \mathcal{E}_h$, we fix a (constant) unit normal vector $\bm{n}_e \in \mathbb{R}^2$ to $e$. In the case where $e \in \mathcal{E}_h^B$, we take $\bm{n}_e$ to be the outward unit normal vector to $\partial\Omega$ restricted to $e$. We denote by $\bm{t}_e$ the unit tangent vector to $e$ obtained by rotating $\bm{n}_e$ counterclockwise by $\pi/2$. The jump operator $\lj \cdot \rj$ across $e$ is then defined as in~\cite{neilan2019discrete}:
\begin{equation*}
\begin{aligned}
\lj v \rj&:= v|_{T_{\text{out}}} - v|_{T_{\text{in}}} \quad &&\text{if} \ e = \partial T_{\text{out}} \cap \partial T_{\text{in}} \in \mathcal{E}_h^I, \\
\lj v \rj&:= v|_{T_{\text{out}}}  &&\text{if} \ e = \partial T_{\text{out}}\cap \partial\Omega  \in \mathcal{E}_h^B,
\end{aligned}
\end{equation*}
where $v$ denotes either a scalar-valued or vector-valued function, and $T_{\text{out}}$ and $T_{\text{in}}$ represent the elements of $\mathcal{T}_h$ sharing the edge $e$, labeled such that $\bm{n}_e$ points outward from $T_{\text{out}}$ and inward toward $T_{\text{in}}$.

Let $\mathcal{P}_k(D)$ denote the space of polynomials of degree at most $k$ on $D$. We define the piecewise polynomial space associated with $\mathcal{T}_h$ by
$$
\mathcal{P}_k(\Omega;\mathcal{T}_h):= \{ v \in L^2(\Omega): v|_{T} \in \mathcal{P}_k(T) \quad \text{for all} \  T \in \mathcal{T}_h \}
$$
and we define the  norm $\|\cdot\|_{L^2(\mathcal{T}_h)}$ by
$$
\|v\|_{L^2(\mathcal{T}_h)}^2 = \sum_{T \in \mathcal{T}_h} \|v\|_{L^2(T)}^2.
$$

The following are definitions of some finite element spaces:
\begin{itemize}
	\item $C^0$ linear finite element space: 
	\begin{align*}
	V_h&:= \{ v_h \in \mathcal{P}_1(\Omega;\mathcal{T}_h) : v_h \ \text{is continuous at all} \ p \in \mathcal{V}_h \}, \\
	V_h^0&:= V_h \cap H_0^1(\Omega), \\
    V_h^g&:= \{ v_h \in V_h: v_h = g_I \ \text{on} \ \partial\Omega \}.
	\end{align*}

	\item Crouzeix--Raviart finite element space:
	$$
	V_h^{CR}:= \{ v_h \in \mathcal{P}_1(\Omega;\mathcal{T}_h) : v_h \ \text{is continuous at all} \ p \in \mathcal{M}_h^{I} \}.
	$$

	\item Morley finite element space:
	\begin{align*}
	M_h:=\Big\{ w_h \in \mathcal{P}_2(\Omega;\mathcal{T}_h): &w_h \ \text{is continuous at all} \  p \in \mathcal{V}_h \ \text{and }  \frac{\partial w_h}{\partial \bm{n}_e} \ \\
	& \text{is continuous at the midpoint of} \ e \ \text{for all} \ e \in \mathcal{E}_h^{I} \Big\}. 
	\end{align*}
\end{itemize}
In addition, we denote by $v_I$ (or $I_hv$) the standard nodal interpolation of the function $v$ onto the $C^0$ linear finite element space, and by $g_I$ the piecewise linear nodal interpolation of $g$ on $\partial\Omega$. Throughout the paper, the symbol $C$ denotes a generic positive constant that is independent of both the mesh size $h$ and the regularization parameter $\epsilon$. We conclude this subsection by listing several auxiliary lemmas that will be used in the subsequent analysis. The proofs of these results can be found in~\cite{brenner2008mathematical}.
\begin{lemma}[Trace and inverse inequalities]\label{lem: trace and inverse inequality}
    For any $T \in \mathcal{T}_h$ and $e \in \mathcal{E}_h$, integers $0 \le n \le m$, and $1 \le p,q \le \infty$, there holds
    \begin{align}
    &\|v\|_{L^2(\partial T)} \le C \Big( h_T^{\frac 12}|v|_{H^1(T)} + h_T^{-\frac12}\|v\|_{L^2(T)} \Big)\quad &&\forall v \in H^1(T),  \label{equ: traceineq}\\
    &\|v\|_{W^{m,q}(T)} \le C h_T^{n-m+ 2\min\{0,1/q-1/p\}}\|v\|_{W^{n,p}(T)} &&\forall v \in \mathcal{P}_k(T), \label{equ: triangle_invineq} \\
    &\|v\|_{W^{m,q}(e)} \le C h_e^{n-m+ \min\{0,1/q-1/p\}}\|v\|_{W^{n,p}(e)} &&\forall v \in \mathcal{P}_k(e). \label{equ: line_invineq}
    \end{align}
\end{lemma}
\begin{lemma}[Bramble--Hilbert Lemma]\label{lem: B-H lem}
    Let $B$ be a ball in $\Omega$ such that $\Omega$ is star-shaped with respect to $B$ and such that its radius $r > \frac12 r_{max}$, where $r_{max}$ is the supremum of the radii of all balls with respect to which $\Omega$ is star-shaped. Let $Q^mu$ be the Taylor polynomial of order $m$ of $u$ averaged over $B$ where $u \in W^{m,p}(\Omega)$ and $p \ge 1$. Then 
    \begin{equation*}
    |u - Q^mu|_{W^{k,p}(\Omega)} \le C d^{m-k}|u|_{W^{m,p}(\Omega)} \qquad k = 0,1,\cdots, m,
    \end{equation*}
    where $d = {\rm diam}(\Omega)$. 
\end{lemma}

\subsection{Recovery operator}

We first recall the construction of the recovery operator $G_h$. For each $v_h \in V_h$, we associate a function $w_h \in M_h$ to $v_h$ as follows:
\begin{itemize}
	\item For all $p \in \mathcal{V}_h$:
    \begin{align*}
        w_h(p) = v_h(p).
    \end{align*}

	\item For all $e \in \mathcal{E}_h$:
    \begin{align*}
        \frac{\partial w_h}{\partial \bm{n}_e}(m_e) = \bm{n}_e \cdot \frac{1}{|\mathcal{K}_e|}\int_{\mathcal{K}_e} \nabla v_h \ \text{d}x,
    \end{align*}
    where $m_e$ is the midpoint of $e$ and $\mathcal{K}_e = \cup_{T \in \mathcal{T}_e}T$.
\end{itemize}
Then, we define the gradient recovery operator $G_h: V_h \rightarrow \mathcal{P}_1(\Omega;\mathcal{T}_h) \times \mathcal{P}_1(\Omega;\mathcal{T}_h)$ by setting $G_hv_h|_{T} = \nabla w_h|_T$ for all $T \in \mathcal{T}_h$.
\begin{lemma}\label{lem: G_h maps V_h to V_h^CR}
	For any $v_h \in V_h$, we have $G_hv_h \in V_h^{CR} \times V_h^{CR}$. Furthermore, $G_hv_h$  has the explicit expression: for all $e \in \mathcal{E}_h$,
	\begin{equation*}
	G_hv_h(m_e) = \frac{1}{|\mathcal{K}_e|} \int_{\mathcal{K}_e}\nabla v_h \ {\rm d}x.
	\end{equation*}
\end{lemma} 

We assume that the recovery operator $G_h$ has the following approximation property.
\begin{assumption}\label{amp: app property of G_h}
	There exists $0 < \gamma < 1$ such that 
	\begin{align*}
		\|\nabla v - G_hv_I\|_{L^2(\Omega)} &\le C h^{1+\gamma} \|v\|_{H^{3}(\Omega)}, \\
		\Big(\sum_{T \in \mathcal{T}_h} h_T^{-2}\|\nabla v - G_hI_hv \|_{L^2(T)}^2 \Big)^{1/2} &\le C h^{\gamma}  \|v\|_{H^3(\Omega)}, \\
		 \|D^2v - DG_hv_I\|_{L^2(\Omega)} &\le C h^{\gamma} \|v\|_{H^{3}(\Omega)},
	\end{align*}
    hold for all $v \in H^{3}(\Omega)$.
\end{assumption}

Assumption \ref{amp: app property of G_h} only depends on the mesh. It holds on most commonly used meshes, see, e.g. \cite[Section 5.1]{chu2025recovery}. We refer Appendix \ref{sec: appendix} and \cite[Section 3.3]{chu2025recovery} for details on the relation between the superconvergence and mesh properties.

\subsection{Some technical results}

In this subsection, we will introduce some results used in the following sections. Let $P_{e}: L^1(e) \to \mathcal{P}_0(e)$ denote the projection defined by 
\begin{equation*}
	P_{e}(v) = \frac{1}{h_e}\int_e v \ \text{d}s.
\end{equation*}  
Then we have the following result.
\begin{lemma}\label{lem: poincare inequality}
	Suppose $\{e_i\}_{i=1}^3$ are the three edges of $T$, then for every $v \in H^1(T)$, we have 
	\begin{equation}
		\sum_{i=1}^3 \| v - P_{e_{i}}(v)\|_{L^2(e_i)} \le C h_T^{1/2} |v|_{H^1(T)}. 
	\end{equation}
\end{lemma}
\begin{proof}
	The result follows from the Bramble--Hilbert lemma on the reference element and a standard scaling argument.
\end{proof}

\begin{lemma}[Lemma 4 and Lemma 5 in \cite{wang2006morley}]\label{lem: edge and triangle}
	For any $v_h \in V_h$, we have 
	\begin{equation*}
		\int_e \lj G_hv_h \rj \ \text{d}s = {\bm 0} \qquad \forall e \in \mathcal{E}_h^I,
	\end{equation*}
	and suppose $e \in \mathcal{E}_h^I$ with $T \cap T^{'} = e$, there holds
	\begin{equation*}
		\big\| \lj G_hv_h \rj \big\|_{L^2(e)} \le C h_T^{1/2} \Big(\big\| D G_hv_h  \big\|_{L^2(T)} + \big\| D G_hv_h  \big\|_{L^2(T^{'})} \Big).
	\end{equation*}
	Moreover if $v_h \in V_h^0$, we have 
	\begin{equation*}
		\big\|  G_hv_h \cdot {\bm t}_e \big\|_{L^2(e)} \le C h_T^{1/2} \big\| D G_hv_h  \big\|_{L^2(T)},
	\end{equation*}
	if $e \in \mathcal{E}_h^B$ with $ e \subset T$.
\end{lemma}

\begin{remark}\label{rem: remark 1}
	A direct corollary of Lemma \ref{lem: edge and triangle} is that, for every $e \in \mathcal{E}_h^I$ and every constant matrix $A \in \R^{2\times 2}$, we have 
	\begin{equation*}
		\int_e \lj G_hv_h \rj \cdot (A{\bm n}_e) \ \text{d}s = 0, \qquad \forall\, v_h \in V_h.
	\end{equation*}
\end{remark}

\begin{lemma}\label{lem: edge}
	Suppose that $A\in \Big(  H^1(\Omega)  \Big)^{2\times 2}$ is a matrix-valued function.
	Then for any $v_h \in V_h$, we have 
	\begin{align*}
		\left| \sum_{e \in \mathcal{E}_h^I}   \int_{e}  \lj G_hv_h  \rj \cdot (A{\bm n}_e) \ \text{d}s    \right|  \le C h \left( \sum_{e\in\mathcal{E}_h^I} h_e^{-1} \big\| \lj G_hv_h \rj \big\|_{L^2(e)}^2    \right)^{1/2} |A |_{H^1(\Omega)}.
	\end{align*}
\end{lemma}
\begin{proof}
	We have the decomposition 
	\begin{align*}
		&\int_{e}  \lj G_hv_h  \rj \cdot (A{\bm n}_e) \ \text{d}s \\
        =&    \int_{e}  \lj G_hv_h  \rj \cdot \Big( \big(A - P_e(A) \big)   {\bm n}_e\Big) \ \text{d}s +    \int_{e}  \lj G_hv_h  \rj \cdot \Big( P_e(A)   {\bm n}_e\Big) \ \text{d}s \\
		 = &\int_{e}  \lj G_hv_h  \rj \cdot \Big( \big(A - P_e(A) \big)   {\bm n}_e\Big) \ \text{d}s,
	\end{align*}
	where we have used Remark \ref{rem: remark 1}. 	Then we have 
	\begin{align*}
		&\Bigg| \sum_{e \in \mathcal{E}_h^I}   \int_{e}  \lj G_hv_h  \rj \cdot (A{\bm n}_e) \ \text{d}s    \Bigg| \\
		 & \le \sum_{e \in \mathcal{E}_h^I} \Bigg|     \int_{e}  \lj G_hv_h  \rj \cdot \Big( \big(A - P_e(A) \big)   {\bm n}_e\Big) \ \text{d}s    \Bigg| \\
		& \le  \sum_{e \in \mathcal{E}_h^I} \|A - P_e(A)\|_{L^2(e)}  \big\| \lj G_hv_h \rj \big\|_{L^2(e)} \\
		& \le  \Big(\sum_{e \in \mathcal{E}_h^I}h_e \|A - P_e(A)\|_{L^2(e)}^2 \Big)^{\frac12}\Big(\sum_{e \in \mathcal{E}_h^I}h_e^{-1}\big\| \lj G_hv_h \rj \big\|_{L^2(e)}^2  \Big)^{\frac12} \\
		& \le Ch \left(\sum_{e \in \mathcal{E}_h^I}h_e^{-1}\big\| \lj G_hv_h \rj \big\|_{L^2(e)}^2  \right)^{1/2}  |A |_{H^1(\Omega)},
	\end{align*}
	where we have used Lemma \ref{lem: poincare inequality} and this completes the  proof.
\end{proof}

\begin{remark}\label{rmk:est}
	For the term $\Big(\sum_{e \in \mathcal{E}_h^I}h_e^{-1}\big\| \lj G_hv_h \rj \big\|_{L^2(e)}^2  \Big)^{1/2}$, by Lemma \ref{lem: edge and triangle} we have 
	\begin{equation*}
		\left(\sum_{e \in \mathcal{E}_h^I}h_e^{-1}\big\| \lj G_hv_h \rj \big\|_{L^2(e)}^2  \right)^{1/2} \le C \|DG_hv_h\|_{L^2(\mathcal{T}_h)}.
	\end{equation*} 
	Moreover, if $v_h = v_I$ for some $v \in H^3(\Omega)$, then by the trace inequality in Lemma \ref{lem: trace and inverse inequality}, 
	\begin{align*}
		&\Big(\sum_{e \in \mathcal{E}_h^I}h_e^{-1}  \big\| \lj G_hv_I \rj \big\|_{L^2(e)}^2  \Big)^{1/2}\\
		  = &\Big(\sum_{e \in \mathcal{E}_h^I}h_e^{-1}\big\| \lj  \nabla v - G_hv_I \rj \big\|_{L^2(e)}^2  \Big)^{1/2} \\
		 \le & C \Big(\sum_{T \in \mathcal{T}_h}    h_T^{-2} \|\nabla v - G_hv_I\|_{L^2(T)}^2 + \sum_{T \in \mathcal{T}_h}  \|D^2v - DG_hv_I\|_{L^2(T)}^2     \Big)^{1/2} 
		 \le  Ch^{\gamma} \|v\|_{H^3(\Omega)}.
	\end{align*}
	The above two estimates will be used in the following sections. Similarly, for $v_h \in V_h^0$ and $v \in H^3(\Omega) \cap H_0^1(\Omega)$ the above estimates are also useful for the term $\Big(\sum_{e \in \mathcal{E}_h^B}h_e^{-1}\big\|  G_hv_h \cdot {\bm t}_e \big\|_{L^2(e)}^2  \Big)^{1/2}$.
\end{remark}

\section{Linearization and finite element approximation}\label{sec:lin}

In this section, we consider the linearization of \eqref{equ: model equations} at the solution $u^{\eps}$ and analyze its finite element approximation. This study lays the groundwork for establishing existence, uniqueness, and error estimates of the numerical scheme for the model problem \eqref{equ: model equations}. We begin by introducing some necessary notation.

We note that for a given smooth function $v$ and $t \in \R$, we have \cite{caffarelli1997properties}
\begin{equation*}
	\det(D^2(u^{\eps} + tv)) = \det(D^2u^{\eps}) + t (\Phi^{\eps}:D^2v) + t^2 \det(D^2v).
\end{equation*}
Thus, the linearization of the operator $M^{\eps}(u^{\eps}):= \eps \Delta^2 u^{\eps} - \det(D^2u^{\eps})$ at $u^{\eps}$ is given by
\begin{equation*}
	L_{u^{\eps}}(v):= \eps\Delta^2v - \Phi^{\eps}:D^2v = \eps\Delta^2v - \text{div}(\Phi^{\eps}\nabla v),
\end{equation*}
where we have used the divergence-free property of the cofactor matrix \cite[Lemma 7]{neilan2010nonconforming}. 

Let $V_0^{*}$ denote the dual of $V_0$. Then for any given $\varphi \in V_0^{*}$, $\psi \in H^{-\frac12}(\partial\Omega)$, we consider the following equations 
\begin{subequations}\label{equ: Linearization of the model equations}
	\begin{align}
		L_{u^{\eps}}(v) = \varphi & \quad \text{in $\Omega$}, \\
		v = 0 & \quad \text{on $\partial\Omega$}, \\
		\Delta v = \psi & \quad \text{on $\partial\Omega$}.
	\end{align}
\end{subequations}

To define the finite element method for \eqref{equ: Linearization of the model equations}, we first consider its weak formulation. Multiplying by a test function $w \in V_0$ and applying the divergence theorem, we obtain the following weak formulation: find $v \in V_0$ such that
\begin{equation}\label{equ: weak formulation of the linearization}
	a(v,w) = \langle\varphi, w \rangle + \eps(\psi, \nabla w\cdot{\bm n})_{L^2(\partial\Omega)} \quad \forall w \in V_0,
\end{equation}
where 
\begin{equation*}
	a(v,w) = \eps(\Delta v, \Delta w)_{L^2(\Omega)} + (\Phi^{\eps} \nabla v, \nabla w)_{L^2(\Omega)}. 
\end{equation*}
Since $u^{\eps}$ is assumed to be convex,  the matrix $\Phi^{\eps}$ is positive semidefinite, which means
\begin{equation*}
	a(v,v) \ge C \eps \|v\|_{H^2(\Omega)}^2. 
\end{equation*}
Hence the well-posedness
 of the problem \eqref{equ: weak formulation of the linearization} follows from the Lax--Milgram Lemma. 

Likewise, to ensure the well-posedness
 of the discrete problem, the following lemma plays a crucial role. Its proof can be found in \cite{chu2025recovery}.
\begin{lemma}[Discrete MT estimate]\label{lem: H^2 bound for V_h^0}
For any $v_h \in V_h^0$, we have the following inequality	
\begin{equation}
\begin{aligned}
	\|DG_hv_h\|_{L^2(\mathcal{T}_h)}^2 \le C\Big( \|{\rm div}(G_hv_h)\|_{L^2(\mathcal{T}_h)}^2 + b_h(v_h,v_h)\Big)
\end{aligned}
\end{equation}
where 
\begin{equation*}
	b_h(u_h,v_h) = \sum_{e \in \mathcal{E}_h^I} h_e^{-1} \big( \lj G_hu_h\rj ,\lj G_hv_h\rj\big)_{L^2(e)} + \sum_{e \in \mathcal{E}_h^B} h_e^{-1} \big( G_hu_h\cdot{\bm t}_e,G_hv_h\cdot {\bm t}_e \big)_{L^2(e)}.
\end{equation*}
\end{lemma}

Inspired by the discrete MT estimate, we design the following finite element method for \eqref{equ: weak formulation of the linearization}: find $v_h \in V_h^0$ such that
\begin{equation}\label{equ: discrete problem for linearization}
	a_h(v_h,w_h) + \eps b_h(v_h,w_h) = \langle\varphi, w_h\rangle+ \eps(\psi, G_hw_h\cdot{\bm n})_{L^2(\partial\Omega)} \quad \forall w_h \in V_h^0,
\end{equation}
where 
\begin{equation*}
	a_h(v_h,w_h) = \eps\big({\rm div}(G_hv_h), {\rm div}(G_hw_h)\big)_{L^2(\Omega)} + (\Phi^{\eps} \nabla v_h, \nabla w_h)_{L^2(\Omega)}. 
\end{equation*}

Then by Lemma \ref{lem: H^2 bound for V_h^0} and the Lax--Milgram lemma we have the following result. 
\begin{theorem}
	Problem \eqref{equ: discrete problem for linearization} has a unique solution $v_h \in V_h^0$.
\end{theorem}

In addition, we can establish the following error estimate:
\begin{theorem}\label{thm: H^2 convergence for linearization}
	Suppose $v \in H^{3}(\Omega)$ satisfies \eqref{equ: Linearization of the model equations}, then we have the following error estimate
	\begin{equation}
		\|D^2v - DG_hv_h\|_{L^2(\mathcal{T}_h)} \le C\Big( h^{\gamma} \|v\|_{H^3(\Omega)} + \eps^{-2} h |v|_{H^2(\Omega)} \Big).
	\end{equation}
\end{theorem}
\begin{proof}
	Let $w_h = v_h - v_I$, then we have 
	\begin{align}
		&a_h(w_h,w_h) + \eps b_h(w_h,w_h) \nonumber\\
        = &a_h(v_h,w_h) + \eps b_h(v_h,w_h) - \big( a_h(v_I,w_h) +\eps b_h(v_I,w_h) \big) \nonumber\\
		= &\langle\varphi, w_h \rangle + \eps(\psi, G_hw_h\cdot{\bm n})_{L^2(\partial\Omega)} \nonumber \\
		& - \big( a_h(v_I,w_h) + \eps b_h(v_I,w_h) \big) \nonumber \\
		= & (\eps\Delta^2v - {\rm div}(\Phi^{\eps}\nabla v) , w_h)_{L^2(\Omega)} + \eps(\psi, G_hw_h\cdot{\bm n})_{L^2(\partial\Omega)} \nonumber \\
		& - \big( a_h(v_I,w_h) +\eps b_h(v_I,w_h) \big). 
	\end{align} 
For the first term, since $w_h \in V_h^0$ we have
\begin{align*}
	&(\eps\Delta^2v - {\rm div}(\Phi^{\eps}\nabla v) , w_h)_{L^2(\Omega)}\\ = & - \eps(\nabla\Delta v, \nabla w_h)_{L^2(\Omega)} + (\Phi^{\eps} \nabla v, \nabla w_h)_{L^2(\Omega)} \\
	= & - \eps(\nabla\Delta v, \nabla w_h - G_hw_h)_{L^2(\Omega)}  - \eps(\nabla\Delta v, G_h w_h)_{L^2(\Omega)} \\
	& + (\Phi^{\eps} \nabla v, \nabla w_h)_{L^2(\Omega)} \\
	= & - \eps(\nabla\Delta v, \nabla w_h - G_hw_h)_{L^2(\Omega)} \\
	& + \eps(\Delta v, {\rm div}(G_hw_h))_{L^2(\Omega)} - \eps(\psi, G_hw_h\cdot{\bm n})_{L^2(\partial\Omega)} \\
	& - \eps\sum_{e \in \mathcal{E}_h^{I}}(\Delta v, \lj G_hw_h\cdot {\bm n}_e \rj )_{L^2(e)} + (\Phi^{\eps} \nabla v, \nabla w_h)_{L^2(\Omega)}.
\end{align*}
This leads to
\begin{align*}
	&a_h(w_h,w_h) +\eps  b_h(w_h,w_h) \\= & \underbrace{\eps(\Delta v - {\rm div}(G_hv_I) , {\rm div}(G_hw_h))_{L^2(\Omega)} + (\Phi^{\eps} \nabla v - \Phi^{\eps} \nabla v_I, \nabla w_h)_{L^2(\Omega)}}_{=:E_1} \\
	& \underbrace{- \eps(\nabla\Delta v, \nabla w_h - G_hw_h)_{L^2(\Omega)}}_{=:E_2} \underbrace{- \eps\sum_{e \in \mathcal{E}_h^{I}}(\Delta v, \lj G_hw_h\cdot {\bm n}_e \rj )_{L^2(e)}}_{=:E_3} \underbrace{- \eps b_h(v_I,w_h)}_{=:E_4}.
\end{align*}
We can deduce the following error bounds 
\begin{align*}
	|E_1| & \le C \eps h^{\gamma} \|v\|_{H^3(\Omega)} \|DG_hw_h\|_{L^2(\mathcal{T}_h)} + C \eps^{-1} h|v|_{H^2(\Omega)} \|\nabla w_h\|_{L^2(\mathcal{T}_h)}, \\
	|E_2| & \le C\eps h|v|_{H^3(\Omega)} \|DG_hw_h\|_{L^2(\mathcal{T}_h)}, \\
	|E_3| & \le C\eps h^{\gamma}|v|_{H^3(\Omega)} \|DG_hw_h\|_{L^2(\mathcal{T}_h)}, \\
	|E_4| & \le C\eps h^{\gamma}\|v\|_{H^3(\Omega)} \|DG_hw_h\|_{L^2(\mathcal{T}_h)}.
\end{align*}
The estimate of $E_3$ follows from Lemma~\ref{lem: edge} and Remark~\ref{rmk:est}, while that of $E_4$ is based on Lemma~\ref{lem: edge and triangle} and Remark~\ref{rmk:est}.
Note that we have 
\begin{equation*}
	a_h(w_h,w_h) + \eps b_h(w_h,w_h) \ge C \eps \|DG_hw_h\|_{L^2(\mathcal{T}_h)}^2. 
\end{equation*}
Then we have the following estimate
\begin{equation}
	\|DG_hw_h\|_{L^2(\mathcal{T}_h)} \le C h^{\gamma} \Big( \|v\|_{H^3(\Omega)} + \eps^{-2} h^{1 - \gamma} |v|_{H^2(\Omega)}  \Big).
\end{equation}
The triangle inequality  implies that 
\begin{equation*}
	\|D^2v - DG_hv_h\|_{L^2(\mathcal{T}_h)} \le \|D^2 v - DG_hv_I\|_{L^2(\mathcal{T}_h)} + \|DG_hw_h\|_{L^2(\mathcal{T}_h)},
\end{equation*}
which concludes our proof. 
\end{proof}

\section{\texorpdfstring{\boldmath $C^0$}{C0} linear FEM for the model problem} \label{sec:nonlinear}

In this section, we present the numerical scheme and its convergence analysis.

\subsection{Numerical scheme}

The weak form of the model equation \eqref{equ: model equations} reads as follows: seek $u^{\eps} \in V_g$ such that for all $v \in V_0$,
\begin{align}
	-\eps(\Delta u^{\eps}, \Delta v)_{L^2(\Omega)} +  \Big( \det(D^2u^{\eps}), v \Big)_{L^2(\Omega)}   = (f,v)_{L^2(\Omega)} - \eps (\eps, \nabla v \cdot {\bm n})_{L^2(\partial\Omega)}.
\end{align}
Then the numerical scheme is designed by the idea of approximating the Hessian by derivatives of the recovered gradient.

Since we consider the nonhomogeneous boundary condition, the definition of $b_h(\cdot,\cdot)$ should be extended to the space $V_h^g \times V_h^0$. We modify the second term in the expression of $b_h$ according to the boundary data $g$. For each $e \in \mathcal{E}_h^B$, let $g_I^2$ denote the quadratic interpolation on edge $e$. Then, for $v_h \in V_h^g$ and $w_h \in V_h^0$, we define
\begin{align*}
    b_h(v_h,w_h) = &\sum_{e \in \mathcal{E}_h^I} h_e^{-1} \big( \lj G_hv_h\rj ,\lj G_hw_h\rj\big)_{L^2(e)} \\
    &+ \sum_{e \in \mathcal{E}_h^B} h_e^{-1} \big( G_hv_h\cdot{\bm t}_e - \partial_{{\bm t}_e}g_I^2,G_hw_h\cdot {\bm t}_e \big)_{L^2(e)}.
\end{align*}
Moreover, if $v_h = v_I$ for some $v \in V_g \cap H^3(\Omega)$, we have 
\begin{equation*}
    |b_h(v_I, w_h)| \le C h^{\gamma} \|v\|_{H^3(\Omega)}\|DG_hw_h\|_{L^2(\mathcal{T}_h)}.
\end{equation*}

Hence we design the following numerical scheme for the model problem: Find $u_h^{\eps} \in V_h^g$ such that for all $v_h \in V_h^0$
\begin{align}\label{equ:newmethod}
	-\eps\Big({\rm div}(G_hu_h^{\eps}), {\rm div}(G_hv_h) \Big)_{L^2(\Omega)} + & \big( \det(DG_hu_h^{\eps}), v_h \big)_{L^2(\Omega)} - \eps b_h (u_h^{\eps}, v_h)  \\
	& = (f,v_h)_{L^2(\Omega)} - \eps (\eps, G_hv_h\cdot {\bm n})_{L^2(\partial\Omega)}. \nonumber
\end{align}
The second term on the left-hand side is nonlinear; hence, the well-posedness of the numerical scheme must be justified, which will be discussed in the next subsection.

\subsection{Existence, uniqueness and $H^2$ error estimate of the numerical scheme}

In this subsection, we prove the well-posedness of the numerical scheme. We first define an operator $T: V_h^{g} \to V_h^{g}$. For any given $v_h \in V_h^{g}$, let $T(v_h) \in V_h^{g}$ denote the solution of the following problem:
\begin{align*}
	&A_h^{\eps}( v_h - T(v_h), w_h ) \\ 
	= & \eps\Big({\rm div}(G_hv_h), {\rm div}(G_hw_h) \Big)_{L^2(\Omega)} - \big( \det(DG_hv_h), w_h \big)_{L^2(\Omega)} \\
	& + \eps b_h (v_h , w_h) + (f,w_h)_{L^2(\Omega)} - (\eps^2, G_hw_h\cdot {\bm n})_{L^2(\partial\Omega)}, \quad \forall\, w_h \in V_h^0,
\end{align*}
where 
\begin{equation*}
	A_h^{\eps}( v_h - T(v_h), w_h ) = a_h(v_h - T(v_h), w_h) + \eps b_h(v_h - T(v_h), w_h).
\end{equation*}

The results of Section~\ref{sec:lin} imply that $T(v_h)$ is well defined. Furthermore, the fixed point of $T$ is the solution of \eqref{equ:newmethod}. To facilitate the error analysis, we introduce the following notation: 
\begin{equation*}
	\mathbb{B}_h(\rho) := \left\{ v_h \in V_h^g: \ \|DG_h(v_h - u_I^{\eps})\|_{L^2(\mathcal{T}_h)} \le \rho \right\}.
\end{equation*}
\begin{lemma}
	We have the following estimates for $u_I^{\eps} - T(u_I^{\eps})$, 
        \begin{align}
		\|DG_h(u_I^{\eps} - T(u_I^{\eps}))\|_{L^2(\mathcal{T}_h)} &\le Ch^{\gamma}\Big( 1 +  \eps^{-1}\|u^{\eps}\|_{W^{2,\infty}(\Omega)}\Big)\|u^{\eps}\|_{H^3(\Omega)} \nonumber  \\
        & \le C_1(\eps^{-1} + \eps^{-3})h^{\gamma},
        \end{align}
	for some constant $C_1 > 0$.
\end{lemma}
\begin{proof}
	Let $w_h^{\eps} = u_I^{\eps} - T(u_I^{\eps})$. Then 
	\begin{align*}
		A_h^{\eps}( w_h^{\eps}, w_h^{\eps} ) = & \eps\Big({\rm div}(G_hu_I^{\eps}), {\rm div}(G_hw_h^{\eps}) \Big)_{L^2(\Omega)} - \big( \det(DG_hu_I^{\eps}), w_h^{\eps} \big)_{L^2(\Omega)} \\
	& + \eps b_h (u_I^{\eps} , w_h^{\eps}) + (f,w_h^{\eps})_{L^2(\Omega)} - (\eps^2, G_hw_h^{\eps}\cdot {\bm n})_{L^2(\partial\Omega)}.
	\end{align*}
	It follows from \eqref{equ: model equations} that 
	\begin{align*}
		&(f,w_h^{\eps})_{L^2(\Omega)}\\
		  = &-\eps(\Delta^2 u^{\eps}, w_h^{\eps})_{L^2(\Omega)} + (\det(D^2u^{\eps}), w_h^{\eps})_{L^2(\Omega)} \\
		 = &\eps (\nabla \Delta u^{\eps}, G_hw_h^{\eps})_{L^2(\Omega)} + \eps (\nabla \Delta u^{\eps}, \nabla w_h^{\eps} - G_hw_h^{\eps})_{L^2(\Omega)} + (\det(D^2u^{\eps}), w_h^{\eps})_{L^2(\Omega)} \\
		 = &-\eps (\Delta u^{\eps}, {\rm div}(G_hw_h^{\eps}))_{L^2(\Omega)} + \eps\sum_{e \in \mathcal{E}_h^{I}}(\Delta u^{\eps}, \lj G_hw_h^{\eps}\cdot {\bm n}_e \rj )_{L^2(e)} + (\eps^2, G_hw_h^{\eps}\cdot {\bm n})_{L^2(\partial\Omega)} \\
		& \quad +\eps (\nabla \Delta u^{\eps}, \nabla w_h^{\eps} - G_hw_h^{\eps})_{L^2(\Omega)} + (\det(D^2u^{\eps}), w_h^{\eps})_{L^2(\Omega)}.
	\end{align*}
	Then we have 
	\begin{align*}
		&A_h^{\eps}( w_h^{\eps}, w_h^{\eps} ) \\
		 = & \eps\Big({\rm div}(G_hu_I^{\eps}) - \Delta u^{\eps}, {\rm div}(G_hw_h^{\eps}) \Big)_{L^2(\Omega)} + \eps\sum_{e \in \mathcal{E}_h^{I}}(\Delta u^{\eps}, \lj G_hw_h^{\eps}\cdot {\bm n}_e \rj )_{L^2(e)} \\
		& \quad  +\eps (\nabla \Delta u^{\eps}, \nabla w_h^{\eps} - G_hw_h^{\eps})_{L^2(\Omega)} + \eps b_h (u_I^{\eps} , w_h^{\eps}) \\
		&\quad  + \big( \det(D^2u^{\eps}) - \det(DG_hu_I^{\eps}), w_h^{\eps} \big)_{L^2(\Omega)}\\
		\le & C \eps h^{\gamma} \|u^{\eps}\|_{H^3(\Omega)}\|DG_hw_h^{\eps}\|_{L^2(\mathcal{T}_h)} + \big( \det(D^2u^{\eps}) - \det(DG_hu_I^{\eps}), w_h^{\eps} \big)_{L^2(\Omega)}.
	\end{align*}
	By the mean value theorem, we have 
	\begin{equation*}
		\big( \det(D^2u^{\eps}) - \det(DG_hu_I^{\eps}), w_h^{\eps} \big)_{L^2(\Omega)} = \Big(\Psi^{\eps}:(D^2u^{\eps} - DG_hu_I^{\eps}), w_h^{\eps}\Big)_{L^2(\Omega)}, 
	\end{equation*}
	where $\Psi^{\eps} = \text{cof}(D^2u^{\eps} - \theta(D^2u^{\eps} - DG_hu_I^{\eps} ))$ for some $\theta \in [0,1]$. Then we have 
	\begin{equation*}
		\|\Psi^{\eps}\|_{L^{\infty}(\Omega)} \le C\|D^2u^{\eps}\|_{L^{\infty}(\Omega)}.
	\end{equation*}
	Using Lemma \ref{lem: discrete Poincare inequality}, we obtain
	\begin{align*}
		A_h^{\eps}( w_h^{\eps}, w_h^{\eps} ) \le  &C \eps h^{\gamma} \|u^{\eps}\|_{H^3(\Omega)}\|DG_hw_h^{\eps}\|_{L^2(\mathcal{T}_h)}\\
		&  + C h^{\gamma} \|u^{\eps}\|_{W^{2,\infty}(\Omega)}\|u^{\eps}\|_{H^3(\Omega)}\|DG_hw_h^{\eps}\|_{L^2(\mathcal{T}_h)}.
	\end{align*}
	Then we can complete the proof by the inequality $A_h^{\eps}( w_h^{\eps}, w_h^{\eps} ) \ge C\eps \|DG_hw_h^{\eps}\|_{L^2(\mathcal{T}_h)}^2$.
\end{proof}

The following lemma gives the contraction estimate for the operator $T$.
\begin{lemma}
	For any $v_h, w_h \in \mathbb{B}_h(\rho)$ there holds
	\begin{equation}
		\|DG_h(T(v_h) - T(w_h))\|_{L^2(\mathcal{T}_h)} \le C(h,\rho,\eps) \|DG_h(v_h - w_h)\|_{L^2(\mathcal{T}_h)},
	\end{equation}
	where $C(h,\rho,\eps) = C_2  \eps^{-1}\Big( \eps^{-1} h^{\gamma}  + \rho \Big)$, for some constant $C_2 > 0$. 
\end{lemma}
\begin{proof}
	Let $z_h = v_h - w_h$. Then by the definition of $T$ we have for any $\eta_h \in V_h^0$:
	\begin{align*}
		&A_h^{\eps}( T(v_h) - T(w_h), \eta_h )\\
		  = & \Big( \det(DG_hv_h) - \det(DG_hw_h), \eta_h   \Big)_{L^2(\Omega)}  + \Big( \Phi^{\eps}\nabla z_h, \nabla \eta_h \Big)_{L^2(\Omega)} \\
		 = & (\Psi_h:DG_hz_h, \eta_h)_{L^2(\Omega)} + \Big( \Phi^{\eps}\nabla z_h, \nabla \eta_h \Big)_{L^2(\Omega)},
	\end{align*}
	where $\Psi_h = \text{cof}(DG_hw_h + \theta(DG_hv_h - DG_hw_h) )$ for some $\theta \in [0,1]$.
	Hence, we have: 
	\begin{align*}
		&A_h^{\eps}( T(v_h) - T(w_h), \eta_h ) \\
		 =  &(\Psi_h:DG_hz_h, \eta_h)_{L^2(\Omega)} + \Big( \Phi^{\eps}G_h z_h, \nabla \eta_h \Big)_{L^2(\Omega)}  + \Big( \Phi^{\eps}(\nabla z_h - G_hz_h), \nabla \eta_h \Big)_{L^2(\Omega)} \\
		 = & \underbrace{\Big( (\Psi_h - \Phi^{\eps}):DG_hz_h, \eta_h \Big)_{L^2(\Omega)}}_{=:R_1} + \underbrace{\sum_{e \in \mathcal{E}_h^I} \Big( \lj \Phi^{\eps} G_hz_h \cdot {\bm n}_e \rj ,\eta_h \Big)_{L^2(e)}}_{=:R_2} \\
		& \quad + \underbrace{\Big( \Phi^{\eps}(\nabla z_h - G_hz_h), \nabla \eta_h \Big)_{L^2(\Omega)}}_{=:R_3}.
	\end{align*}
	Then we have 
	\begin{align*}
	|R_2| &\le  C  h \|DG_hz_h\|_{L^2(\mathcal{T}_h)}|\Phi^{\eps} \eta_h  |_{H^1(\Omega)} \\
	&\le Ch  \|DG_hz_h\|_{L^2(\mathcal{T}_h)} \Big(  |\Phi^{\eps}   |_{H^1(\Omega)} \|\eta_h\|_{L^{\infty}(\Omega)} + |\eta_h|_{H^1(\Omega)}\|\Phi^{\eps}\|_{L^{\infty}(\Omega)}        \Big)          \\
	&\le C(\|u^{\eps}\|_{W^{2,\infty}(\Omega)}+ \|u^{\eps}\|_{H^3(\Omega)}) h \|DG_hz_h\|_{L^2(\mathcal{T}_h)} \|DG_h \eta_h\|_{L^2(\mathcal{T}_h)},     
    \end{align*}
    and
    \begin{align*}
		|R_3| &\le C \|u^{\eps}\|_{W^{2,\infty}(\Omega)} h \|DG_hz_h\|_{L^2(\mathcal{T}_h)} \|DG_h \eta_h\|_{L^2(\mathcal{T}_h)}.
	\end{align*}
	Moreover, we have 
	\begin{align*}
		&\|\Psi_h - \Phi^{\eps}\|_{L^2(\Omega)} \\
		 \le & \|D^2u^{\eps} - DG_hu_I^{\eps}\|_{L^2(\Omega)} + \|DG_hu_I^{\eps} - DG_hv_h\|_{L^2(\Omega)} + \|DG_hz_h\|_{L^2(\Omega)} \\
		 \le &C (h^{\gamma} \|u^{\eps}\|_{H^3(\Omega)} + \rho ).
	\end{align*}
	Hence, we have 
	\begin{align*}
		|R_1| & \le C (h^{\gamma} \|u^{\eps}\|_{H^3(\Omega)} + \rho ) \|DG_hz_h\|_{L^2(\mathcal{T}_h)} \|\eta_h \|_{L^{\infty}(\Omega)} \\
		&\le C  (h^{\gamma} \|u^{\eps}\|_{H^3(\Omega)} + \rho ) \|DG_hz_h\|_{L^2(\mathcal{T}_h)} \|DG_h\eta_h\|_{L^2(\mathcal{T}_h)}.
	\end{align*}
	Then we complete the proof by the above estimates.
\end{proof}
\begin{theorem}
	For $h \le h_0$, the numerical scheme \eqref{equ:newmethod} has a unique solution $u_h^{\eps} \in \mathbb{B}_h(\rho)$ and we have the following $H^2$ error estimate 
        \begin{align}
		\|D^2 u^{\eps} - DG_h u_h^{\eps}\|_{L^2(\mathcal{T}_h)} &\le Ch^{\gamma}\Big( 1 +  \eps^{-1}\|u^{\eps}\|_{W^{2,\infty}(\Omega)}\Big)\|u^{\eps}\|_{H^3(\Omega)} \nonumber \\
        & \le C \Big( \eps^{-1} +  \eps^{-3}\Big) h^{\gamma}.
        \end{align}
\end{theorem}
\begin{proof}
	Set $\rho_0 = 2C_1h^{\gamma}\Big( 1 +  \eps^{-1}\|u^{\eps}\|_{W^{2,\infty}(\Omega)}\Big)\|u^{\eps}\|_{H^3(\Omega)}$, and take $h_0$ sufficiently small such that 
	$C(h, \rho_0, \eps) \le 1/2$ when $h \le h_0$. Then we have 
	\begin{align*}
		\|DG_h(T(v_h) - T(w_h))\|_{L^2(\mathcal{T}_h)} \le \frac12 \|DG_h(v_h - w_h)\|_{L^2(\mathcal{T}_h)} \quad \forall v_h,w_h \in \mathbb{B}_h(\rho_0),
	\end{align*}
	and for any $z_h \in \mathbb{B}_h(\rho_0)$, 
	\begin{align*}
		&\|DG_hT(z_h) - DG_hu_I^{\eps}\|_{L^2(\mathcal{T}_h)}\\
		  \le &\|DG_hT(z_h) - DG_hT(u_I^{\eps})\|_{L^2(\mathcal{T}_h)} + \|DG_hT(u_I^{\eps}) - DG_hu_I^{\eps}\|_{L^2(\mathcal{T}_h)} \\
		 \le &\frac12 \rho_0 + \frac12 \rho_0 = \rho_0.
	\end{align*}
	Hence, by the Banach fixed-point theorem, there exists a unique $u_h^{\eps} \in \mathbb{B}_h(\rho_0)$ such that $T(u_h^{\eps}) = u_h^{\eps}$ which means $u_h^{\eps}$ is a solution of the numerical scheme. Meanwhile we have 
	\begin{align*}
		\|D^2u^{\eps} - DG_hu_h^{\eps}\|_{L^2(\mathcal{T}_h)} & \le \|D^2u^{\eps} - DG_hu_I^{\eps}\|_{L^2(\mathcal{T}_h)} + \|DG_hu_I^{\eps} - DG_hu_h^{\eps}\|_{L^2(\mathcal{T}_h)} \\
		& \le Ch^{\gamma} \|u^{\eps}\|_{H^3(\Omega)} + \rho_0 \\
		& \le Ch^{\gamma}\Big( 1 +  \eps^{-1}\|u^{\eps}\|_{W^{2,\infty}(\Omega)}\Big)\|u^{\eps}\|_{H^3(\Omega)},
	\end{align*}
	which completes the proof.
\end{proof}

\subsection{\texorpdfstring{\boldmath $H^1$}{H1} error estimate}

In this subsection, we give an $H^1$ error estimate by the duality argument. It is well known that the duality argument depends on the regularity of the solution of the auxiliary problem:
\begin{subequations}\label{equ: Duality argument_1}
	\begin{align}
		L_{u^{\eps}}(\mathcal{U}) = \varphi & \quad \text{in $\Omega$}, \\
		\mathcal{U} = 0 & \quad \text{on $\partial\Omega$}, \\
		\Delta \mathcal{U} = 0 & \quad \text{on $\partial\Omega$},
	\end{align}
\end{subequations}
where $\varphi \in H^{-1}(\Omega) \subset V_0^{*}$ is defined as 
\begin{equation*}
	\langle \varphi, w \rangle = \Big(\nabla (u_h^{\eps} - u_I^{\eps}), \nabla w \Big)_{L^2(\Omega)} \quad \forall w \in H_0^{1}(\Omega).
\end{equation*}
As in \cite{neilan2010nonconforming}, we can make the following assumption about the auxiliary problem. 
\begin{assumption}
	We assume that the solution of \eqref{equ: Duality argument_1} satisfies the following regularity property: 
	\begin{equation*}
		\|\mathcal{U}\|_{H^3(\Omega)} \le C_s(\eps) \|\varphi\|_{H^{-1}(\Omega)} \le C_s(\eps) \|\nabla (u_h^{\eps} - u_I^{\eps})\|_{L^2(\Omega)}.
	\end{equation*}
\end{assumption}

The weak form of \eqref{equ: Duality argument_1} reads as follows: 
 find $\mathcal{U} \in V_0$ such that 
\begin{equation*}
	a(w, \mathcal{U}) = (\nabla (u_h^{\eps} - u_I^{\eps}), \nabla w)_{L^2(\Omega)} \quad \forall w \in V_0.
\end{equation*}
The discrete auxiliary problem is to find $\mathcal{U}_h \in V_h^0$ such that
\begin{equation*}
	A_h^{\eps}(w_h, \mathcal{U}_h) = (\nabla (u_h^{\eps} - u_I^{\eps}), \nabla w_h)_{L^2(\Omega)} \quad \forall w_h \in V_h^0. 
\end{equation*}
Thus, we have 
\begin{equation*}
    \begin{aligned}
     \|\mathcal{U}\|_{H^2(\Omega)} &\le C \eps^{-1} \|\nabla (u_h^{\eps} - u_I^{\eps})\|_{L^2(\Omega)}, \\
	 \|DG_h\mathcal{U}_h\|_{L^2(\mathcal{T}_h)} &\le C \eps^{-1} \|\nabla (u_h^{\eps} - u_I^{\eps})\|_{L^2(\Omega)},
    \end{aligned}
\end{equation*}
and by Theorem \ref{thm: H^2 convergence for linearization}
\begin{equation*}
	\|D^2\mathcal{U} - DG_h\mathcal{U}_h\|_{L^2(\mathcal{T}_h)} \le C \Big(C_s(\eps) h^{\gamma} + \eps^{-3}h \Big) \|\nabla (u_h^{\eps} - u_I^{\eps})\|_{L^2(\Omega)}.
\end{equation*}

\begin{theorem}\label{thm: H^1 error}
	The solution $u_h^{\eps}$ satisfies the following error estimate
        \begin{align}
            \|\nabla (u_h^{\eps} - u_I^{\eps})\|_{L^2(\Omega)} & \le  C(\eps) ( h^{2\gamma} + h^{1 + \gamma} + h^2 )    + C \eps^{-1} h ,\label{equ: H^1 error}
        \end{align}
	and this implies the following $H^1$ error estimate 
        \begin{align}
            \|\nabla ( u^{\eps} - u_h^{\eps} )\|_{L^2(\Omega)} & \le C(\eps) ( h^{2\gamma} + h^{1 + \gamma} + h^2 )   + C \Big( \eps^{-\frac12} +  \eps^{-1} \Big) h.
        \end{align}
\end{theorem}
\begin{proof}
	We have the decomposition
	\begin{align*}
		\|\nabla (u_h^{\eps} - u_I^{\eps}) \|_{L^2(\Omega)}^2  = A_h^{\eps} (u_h^{\eps} - u_I^{\eps}, \mathcal{U}_h) = A_h^{\eps} (u_h^{\eps} , \mathcal{U}_h) - A_h^{\eps} (u_I^{\eps}, \mathcal{U}_h)  =: \sum_{i = 1}^6 \mathcal{G}_i,
	\end{align*}
	where 
	\begin{align*}
		\mathcal{G}_1 & = \eps ( \Delta u^{\eps} - {\rm div}(G_hu_I^{\eps}),   {\rm div}(G_h\mathcal{U}_h)  )_{L^2(\Omega)} \\
		\mathcal{G}_2 &= - \eps(\nabla\Delta u^{\eps}, \nabla\mathcal{U}_h - G_h\mathcal{U}_h)_{L^2(\Omega)} \\
		\mathcal{G}_3 & = - \eps\sum_{e \in \mathcal{E}_h^{I}}(\Delta u^{\eps}, \lj G_h\mathcal{U}_h\cdot {\bm n}_e \rj )_{L^2(e)} \\
		\mathcal{G}_4 & = (\Phi^{\eps} \nabla (u^{\eps} -  u_I^{\eps}), \nabla \mathcal{U}_h)_{L^2(\Omega)} \\
		\mathcal{G}_5 & = - \Big( 
		\det(D^2u^{\eps}) - \det(DG_hu^{\eps}_h), \mathcal{U}_h \Big)_{L^2(\Omega)} + (\Phi^{\eps} \nabla (u^{\eps}_h -  u^{\eps}), \nabla \mathcal{U}_h)_{L^2(\Omega)} \\
		\mathcal{G}_6 & = -\eps b_h(u_I^{\eps}, \mathcal{U}_h).
	\end{align*}
	Then we can obtain the following error bounds 
	\begin{align*}
		|\mathcal{G}_2|  \le  C  h \|u^{\eps}\|_{H^3(\Omega)}\|\nabla (u_h^{\eps} - u_I^{\eps})\|_{L^2(\Omega)} \le C \eps^{-1} h \|\nabla (u_h^{\eps} - u_I^{\eps})\|_{L^2(\Omega)}, 
    \end{align*}
    and
    \begin{align*}
		|\mathcal{G}_3| \le  C\eps h^{1 + \gamma} \|u^{\eps}\|_{H^3(\Omega)}\|\mathcal{U}\|_{H^3(\Omega)}   \le CC_s(\eps)  h^{1 + \gamma} \|\nabla (u_h^{\eps} - u_I^{\eps})\|_{L^2(\Omega)}, 
    \end{align*}
    and
    \begin{align*}
		|\mathcal{G}_6| \le & \eps |b_h(u_I^{\eps} , \mathcal{U}_I)| + \eps| b_h(u_I^{\eps} , \mathcal{U}_h - \mathcal{U}_I) |      \\ 
		 \le & C\eps h^{2\gamma} \|u^{\eps}\|_{H^3(\Omega)}\|\mathcal{U}\|_{H^3(\Omega)} + C\eps h^{\gamma} \|u^{\eps}\|_{H^3(\Omega)} \|DG_h(\mathcal{U}_h - \mathcal{U}_I) \|_{L^2(\mathcal{T}_h)} \\
        \le &C C_s(\eps) h^{2\gamma} \|\nabla (u_h^{\eps} - u_I^{\eps})\|_{L^2(\Omega)}  \\
        &+ C h^{\gamma} \Big(  \|D^2 \mathcal{U} - DG_h\mathcal{U}_h  \|_{L^2(\mathcal{T}_h)} +  \|D^2 \mathcal{U} - DG_h\mathcal{U}_I \|_{L^2(\mathcal{T}_h)}  \Big)\\
		\le & C \Big( C_s(\eps) h^{2\gamma} + \eps^{-3} h^{1 + \gamma}\Big) \|\nabla (u_h^{\eps} - u_I^{\eps})\|_{L^2(\Omega)}.
	\end{align*}
	
	For $\mathcal{G}_1$, we have 
	\begin{align*}
		\mathcal{G}_1 & = \eps ( \Delta u^{\eps} - {\rm div}(G_hu_I^{\eps}),   {\rm div}(G_h\mathcal{U}_h) - \Delta \mathcal{U}  )_{L^2(\Omega)} + \eps ( \Delta u^{\eps} - {\rm div}(G_hu_I^{\eps}),   \Delta \mathcal{U}  )_{L^2(\Omega)} \\
		&= \eps ( \Delta u^{\eps} - {\rm div}(G_hu_I^{\eps}),   {\rm div}(G_h\mathcal{U}_h) - \Delta \mathcal{U}  )_{L^2(\Omega)} - \eps (\nabla u^{\eps} - G_hu_I^{\eps}, \nabla \Delta \mathcal{U}   )_{L^2(\Omega)} \\
		&\quad + \eps \sum_{e \in \mathcal{E}_h^I} (\lj (\nabla u^{\eps} - G_h u_I^{\eps})\cdot{\bm n}_e \rj, \Delta \mathcal{U})_{L^2(e)}.
	\end{align*}
	Hence, we have 
	\begin{align*}
		| \mathcal{G}_1 | & \le C \eps h^{2\gamma} \|u^{\eps}\|_{H^3(\Omega)}\|\mathcal{U}\|_{H^3(\Omega)} + C\eps h^{1 + \gamma} \|u^{\eps}\|_{H^3(\Omega)}\|\mathcal{U}\|_{H^3(\Omega)} \\
		&\le CC_s(\eps)  h^{2\gamma} \|\nabla (u_h^{\eps} - u_I^{\eps})\|_{L^2(\Omega)}.
	\end{align*}
	
	For $\mathcal{G}_4$, we have
	\begin{align*}
		\mathcal{G}_4 & = (\Phi^{\eps} \nabla (u^{\eps} -  u_I^{\eps}), \nabla \mathcal{U}_h)_{L^2(\Omega)} = (\nabla (u^{\eps} -  u_I^{\eps}), \Phi^{\eps} \nabla \mathcal{U}_h)_{L^2(\Omega)} \\
		& = (\nabla (u^{\eps} -  u_I^{\eps}), \Phi^{\eps} G_h\mathcal{ U}_h)_{L^2(\Omega)} + (\nabla (u^{\eps} -  u_I^{\eps}), \Phi^{\eps} \nabla \mathcal{U}_h - \Phi^{\eps} G_h \mathcal{U}_h)_{L^2(\Omega)} \\
		& = \underbrace{(\nabla (u^{\eps} -  u_I^{\eps}), \Phi^{\eps} \nabla \mathcal{U}_h - \Phi^{\eps} G_h \mathcal{U}_h)_{L^2(\Omega)} - (u^{\eps} -  u_I^{\eps}, \Phi^{\eps}:D G_h \mathcal{U}_h)_{L^2(\Omega)}}_{=:\mathcal{G}_{41}} \\
		&\quad + \underbrace{\sum_{e \in \mathcal{E}_h^I} (\lj \Phi^{\eps} G_h \mathcal{U}_h\cdot{\bm n}_e \rj, u^{\eps} -  u_I^{\eps})_{L^2(e)}}_{=:\mathcal{G}_{42}} + \underbrace{\sum_{e \in \mathcal{E}_h^{B}} ( g - g_I, \Phi^{\eps} G_h \mathcal{U}_h\cdot{\bm n}_e)_{L^2(e)}}_{=:\mathcal{G}_{43}} .
	\end{align*}
	We can estimate $\mathcal{G}_{41}$ as
	\begin{align*}
		|\mathcal{G}_{41}|  & \le C \|u^{\eps}\|_{W^{2,\infty}(\Omega)} h^2 \|u^{\eps}\|_{H^2(\Omega)}\|DG_h\mathcal{U}_h\|_{L^2(\mathcal{T}_h)}\\
		& \le C \eps^{-\frac52} h^2\|\nabla (u_h^{\eps} - u_I^{\eps})\|_{L^2(\Omega)}.
	\end{align*}
	In addition, we have
	\begin{align*}
		|\mathcal{G}_{42}| & \le C h \|DG_h\mathcal{U}_h\|_{L^2(\mathcal{T}_h)} |(u^{\eps} - u_I^{\eps}) \Phi^{\eps}|_{H^1(\Omega)} \\
		& \le C h \|DG_h\mathcal{U}_h\|_{L^2(\mathcal{T}_h)} \Big( \| u^{\eps} - u_I^{\eps}\|_{L^{\infty}(\Omega)}| \Phi^{\eps}|_{H^1(\Omega)} + \|  \Phi^{\eps} \|_{L^{\infty}(\Omega)} |u^{\eps} - u_I^{\eps}|_{H^1(\Omega)}   \Big) \\
		& \le C \Big(  \eps^{-2} + \eps^{-\frac52}   \Big) h^2\|\nabla (u_h^{\eps} - u_I^{\eps})\|_{L^2(\Omega)},
	\end{align*}
    and 
    \begin{align*}
        |\mathcal{G}_{43}| & \le C\|\Phi^{\eps}\|_{L^{\infty}(\Omega)} \|g - g_I\|_{L^2(\partial\Omega)}\|DG_h\mathcal{U}_h\|_{L^2(\mathcal{T}_h)} \\
        & \le C\eps^{-3} h^2 \|\nabla (u_h^{\eps} - u_I^{\eps})\|_{L^2(\Omega)}.
    \end{align*}
	These estimates imply 
	\begin{equation*}
		|\mathcal{G}_{4}|  \le C \Big(  \eps^{-2} + \eps^{-\frac52} + \eps^{-3}   \Big) h^2\|\nabla (u_h^{\eps} - u_I^{\eps})\|_{L^2(\Omega)}.
	\end{equation*}

	For $\mathcal{G}_5$, we have
    {\small
	\begin{align*}
		 \mathcal{G}_5 & = - \Big( 
		\det(D^2u^{\eps}) - \det(DG_hu^{\eps}_h), \mathcal{U}_h \Big)_{L^2(\Omega)} + (\Phi^{\eps} (\nabla u^{\eps}_h -  G_h u_h^{\eps}), \nabla \mathcal{U} _h)_{L^2(\Omega)} \\
		& \quad + (\Phi^{\eps} (G_hu_h^{\eps} -\nabla u^{\eps}), \nabla \mathcal{U}_h)_{L^2(\Omega)} \\
		& = - \Big(\det(D^2u^{\eps}) - \det(DG_hu^{\eps}_h), \mathcal{U}_h \Big)_{L^2(\Omega)} + (\Phi^{\eps} (\nabla u^{\eps}_h -  G_h u_h^{\eps}), \nabla \mathcal{U} _h)_{L^2(\Omega)} \\
		&\quad - (\Phi^{\eps}: (DG_hu_h^{\eps} -D^2u^{\eps}), \mathcal{U}_h)_{L^2(\Omega)} - \sum_{e \in \mathcal{E}_h^I} (\lj \Phi^{\eps}(\nabla u^{\eps} - G_h u_h^{\eps})\cdot{\bm n}_e \rj, \mathcal{U}_h)_{L^2(e)} \\
		& = \Big( (\Psi_h - \Phi^{\eps}):(DG_hu_h^{\eps} -D^2u^{\eps}), \mathcal{U}_h  \Big)_{L^2(\Omega)} + (\Phi^{\eps} (\nabla u^{\eps}_h -  G_h u_h^{\eps}), \nabla \mathcal{U} _h)_{L^2(\Omega)} \\
		&\quad - \sum_{e \in \mathcal{E}_h^I} (\lj \Phi^{\eps}(\nabla u^{\eps} - G_h u_h^{\eps})\cdot{\bm n}_e \rj, \mathcal{U}_h)_{L^2(e)} \\
		& = \underbrace{\Big( (\Psi_h - \Phi^{\eps}):(DG_hu_h^{\eps} -D^2u^{\eps}), \mathcal{U}_h  \Big)_{L^2(\Omega)}}_{=:\mathcal{G}_{51}}  \underbrace{- \sum_{e \in \mathcal{E}_h^I} (\lj \Phi^{\eps}(\nabla u^{\eps} - G_h u_h^{\eps})\cdot{\bm n}_e \rj, \mathcal{U}_h)_{L^2(e)}}_{=:\mathcal{G}_{52}} \\
		&\quad + \underbrace{(\Phi^{\eps} (\nabla u^{\eps}_h -  G_h u_h^{\eps}), \nabla \mathcal{U}_h - \nabla \mathcal{U})_{L^2(\Omega)}}_{=:\mathcal{G}_{53}} + \underbrace{(\Phi^{\eps} (\nabla u^{\eps}_h -  G_h u_h^{\eps}), \nabla \mathcal{U} )_{L^2(\Omega)}}_{=:\mathcal{G}_{54}},
	\end{align*}}
where $\Psi_h = \text{cof}(D^2u^{\eps} + \theta (DG_hu_h^{\eps} - D^2u^{\eps}))$ for some $\theta$ in $[0,1]$. We estimate the four terms one by one. Firstly, we have
	\begin{align*}
		|\mathcal{G}_{51} | & \le C \|D^2u^{\eps} - DG_hu_h^{\eps}\|_{L^2(\Omega)}^2 \|\mathcal{U}_h\|_{L^{\infty}(\Omega)} \\
		&\le C   \|D^2u^{\eps} - DG_hu_h^{\eps}\|_{L^2(\Omega)}^2 \|DG_h\mathcal{U}_h\|_{L^2(\T_h)} \\
		& \le C  \Big( \eps^{-3} +  \eps^{-5} + \eps^{-7}\Big) h^{2\gamma} \|\nabla (u_h^{\eps} - u_I^{\eps})\|_{L^2(\Omega)},
	\end{align*}
	and 
	\begin{align*}
		|\mathcal{G}_{52} | & \le C h \|DG_h(u_h^{\eps} - u_I^{\eps}) \|_{L^2(\mathcal{T}_h)} |\Phi^{\eps}\mathcal{U}_h|_{H^1(\Omega)} + Ch^{1 + \gamma} \|u^{\eps}\|_{H^3(\Omega)}|\Phi^{\eps}\mathcal{U}_h|_{H^1(\Omega)} \\
		& \le C \Big( h^{1 + \gamma} \|u^{\eps}\|_{H^3(\Omega)} + h\|D^2u^{\eps} - DG_hu_h^{\eps}\|_{L^2(\Omega)} \Big) \\
        & \quad \times \Big(   \| \mathcal{U}_h\|_{L^{\infty}(\Omega)}| \Phi^{\eps}|_{H^1(\Omega)} + \|  \Phi^{\eps} \|_{L^{\infty}(\Omega)} |\mathcal{U}_h|_{H^1(\Omega)}   \Big) \\
		& \le C \eps^{-2} \Big( \eps^{-1} h^{1 + \gamma} + h^{1 + \gamma} (\eps^{-1} + \eps^{-3}) \Big) \|\nabla (u_h^{\eps} - u_I^{\eps})\|_{L^2(\Omega)} \\
        & \le C \Big(  \eps^{-3} + \eps^{-5}  \Big)  h^{1 + \gamma} \|\nabla (u_h^{\eps} - u_I^{\eps})\|_{L^2(\Omega)}.
	\end{align*}
Then, we estimate $\mathcal{G}_{53}$ as
	\begin{align*}
		|\mathcal{G}_{53} | & \le C \|u^{\eps}\|_{W^{2,\infty}(\Omega)} h \|DG_hu_h^{\eps}\|_{L^2(\mathcal{T}_h)} \|DG_h(\mathcal{U}_h - \mathcal{U}_I)\|_{L^2(\mathcal{T}_h)} \\
        & \quad + C\|u^{\eps}\|_{W^{2,\infty}(\Omega)}h^2 \|DG_hu_h^{\eps}\|_{L^2(\mathcal{T}_h)} |\mathcal{U}|_{H^2(\Omega)} \\ 
        & \le C \eps^{-1} h \|DG_hu_h^{\eps}\|_{L^2(\mathcal{T}_h)} \Big(   \|D^2\mathcal{U} - DG_h\mathcal{U}_h \|_{L^2(\mathcal{T}_h)} + \|D^2\mathcal{U} - DG_h\mathcal{U}_I\|_{L^2(\mathcal{T}_h)} \Big) \\ 
        & \quad + C \eps^{-2} h^2 \|DG_hu_h^{\eps}\|_{L^2(\mathcal{T}_h)} \|\nabla (u_h^{\eps} - u_I^{\eps})\|_{L^2(\Omega)} \\
        & \le C \eps^{-1} h  \Big(C_s(\eps) h^{\gamma} + \eps^{-3}h \Big)  \|DG_hu_h^{\eps}\|_{L^2(\mathcal{T}_h)} \|\nabla (u_h^{\eps} - u_I^{\eps})\|_{L^2(\Omega)} \\
        & \quad + C \eps^{-2} h^2 \|DG_hu_h^{\eps}\|_{L^2(\mathcal{T}_h)} \|\nabla (u_h^{\eps} - u_I^{\eps})\|_{L^2(\Omega)} \\
        & \le C \Big(     \eps^{-1 } C_s(\eps) + \eps^{-2} + \eps^{-4}   \Big) h^{1 + \gamma} \|DG_hu_h^{\eps}\|_{L^2(\mathcal{T}_h)} \|\nabla (u_h^{\eps} - u_I^{\eps})\|_{L^2(\Omega)}.
	\end{align*}
Recalling the $H^2$ error estimate
    \begin{equation*}
		\|D^2 u^{\eps} - DG_h u_h^{\eps}\|_{L^2(\mathcal{T}_h)} \le C \Big( \eps^{-1} +  \eps^{-3}\Big) h^{\gamma},
	\end{equation*}
we obtain, for $h$ sufficiently small,
    \begin{equation*}
        \|DG_h u_h^{\eps}\|_{L^2(\mathcal{T}_h)} \le C \eps^{-\frac12}, 
    \end{equation*}
which implies that
    \begin{equation*}
        |\mathcal{G}_{53} | \le  C \Big(     \eps^{-\frac32 } C_s(\eps) + \eps^{-\frac52} + \eps^{-\frac92}   \Big) h^{1 + \gamma} \|\nabla (u_h^{\eps} - u_I^{\eps})\|_{L^2(\Omega)}.
    \end{equation*}
Moreover,
	\begin{align*}
		|\mathcal{G}_{54} | & \le C h^2 \|DG_hu_h^{\eps}\|_{L^2(\mathcal{T}_h)} |\Phi^{\eps} \nabla \mathcal{U}|_{H^1(\Omega)} \\
		& \le C \eps^{-\frac12} h^2   \Big(  \|\Phi^{\eps}\|_{L^{\infty}(\Omega)} |\mathcal{U}|_{H^2(\Omega)} + |\Phi^{\eps}|_{H^1(\Omega)} \|\nabla \mathcal{U}\|_{L^{\infty}(\Omega)} \Big) \\
		& \le C \Big(   \eps^{-\frac52} + \eps^{-\frac32} C_s(\eps)    \Big) h^2 \|\nabla (u_h^{\eps} - u_I^{\eps})\|_{L^2(\Omega)}.
	\end{align*}
The above estimates imply that
    \begin{align*}
		|\mathcal{G}_5|  & \le  C  \Big(   \eps^{-3} +   \eps^{-5} + \eps^{-7}    \Big)  h^{2\gamma} \|\nabla (u_h^{\eps} - u_I^{\eps})\|_{L^2(\Omega)} \\
        & \quad +  C  \Big(   \eps^{-3} +   \eps^{-5} + \eps^{-\frac32 } C_s(\eps) + \eps^{-\frac52} + \eps^{-\frac92}    \Big) h^{1 + \gamma} \|\nabla (u_h^{\eps} - u_I^{\eps})\|_{L^2(\Omega)} \\
        & \quad + C \Big(   \eps^{-\frac52} + \eps^{-\frac32} C_s(\eps)    \Big) h^2 \|\nabla (u_h^{\eps} - u_I^{\eps})\|_{L^2(\Omega)}.
	\end{align*}
Consequently, we have
	\begin{align*}
		\|\nabla (u_h^{\eps} - u_I^{\eps})\|_{L^2(\Omega)} \le  C(\eps) ( h^{2\gamma} + h^{1 + \gamma} + h^2 )    + C \eps^{-1} h,
	\end{align*}
and the proof is completed by applying the triangle inequality.
\end{proof}
\begin{remark}
	The estimate in \eqref{equ: H^1 error} is not sharp if $\gamma > 1/2$ and $u^{\eps}$ satisfies $\Delta u^{\eps} \in H^2(\Omega)$, due to the estimate of the term $\mathcal{G}_2$ in Theorem~\ref{thm: H^1 error}. Under the assumption that $u^{\eps} \in H^3(\Omega)$, we expect $\Delta u^{\eps} \in H^2(\Omega)$, since $\Delta u^{\eps}$ can be viewed as the solution of the following equation: 
	\begin{equation*}
		- \epsilon \Delta z =  f - \det(D^2u^{\eps}) \quad \text{in $\Omega$},\quad  z = \eps  \quad \text{on $\partial\Omega$}.
	\end{equation*}
Hence, the estimate of $\mathcal{G}_2$ can be improved by Theorem \ref{thm: superconvergence property of G_h} as
	\begin{align*}
		|\mathcal{G}_2| &\le C |\Delta u^{\eps}|_{H^2(\Omega)} h^2 \|\nabla (u_h^{\eps} - u_I^{\eps})\|_{L^2(\Omega)} \\
        & \le C \eps^{-1} \Big(  \|f\|_{L^2(\Omega)} + \|\det(D^2u^{\eps})\|_{L^2(\Omega)}  \Big)h^2 \|\nabla (u_h^{\eps} - u_I^{\eps})\|_{L^2(\Omega)} \\
        & \le C \eps^{-3} h^2 \|\nabla (u_h^{\eps} - u_I^{\eps})\|_{L^2(\Omega)}.
	\end{align*}
Then, the estimate of $\|\nabla (u_h^{\eps} - u_I^{\eps})\|_{L^2(\Omega)}$ can be improved as
\begin{equation*}
		\|\nabla (u_h^{\eps} - u_I^{\eps})\|_{L^2(\Omega)} \le C(\eps) (h^{2\gamma} + h^{1 + \gamma} + h^2),
	\end{equation*}
and this implies the following estimate:
	\begin{equation*}
		\|\nabla u^{\eps} - G_hu_h^{\eps}\|_{L^2(\Omega)} \le C \eps^{-1} h^{1 + \gamma} +  C(\eps)(h^{2\gamma} + h^{1 + \gamma} + h^2).
	\end{equation*}

    By the Poincar\'e inequality, we have 
    \begin{align*}
        \|u_h^{\eps} - u_I^{\eps}\|_{L^2(\Omega)} \le \|\nabla (u_h^{\eps} - u_I^{\eps})\|_{L^2(\Omega)} \le C(\eps) (h^{2\gamma} + h^{1 + \gamma} + h^2).
    \end{align*}
    This implies an $L^2$ estimate
    \begin{align*}
        \|u^{\eps} - u_h^{\eps}\|_{L^2(\Omega)} \le C \eps^{-\frac12}h^2 + C(\eps) (h^{2\gamma} + h^{1 + \gamma} + h^2).
    \end{align*}
\end{remark}

\section{Numerical experiments}\label{sec:num}
 
In this section, we provide several experiments to gauge the efficiency of the finite element method developed in the previous sections. Unless otherwise stated, all computations are performed on the domain $\Omega = (0,1)^2$, and Newton’s method is used as the nonlinear solver. The numerical tests are carried out on both uniform and unstructured meshes. The numerical error will be quantified in the following norms:
\begin{equation*}
	e = \|u_h - u\|_{L^2(\Omega)}, \quad De = \|\nabla u_h - \nabla u\|_{L^2(\mathcal{T}_h)},\quad D^2e = \|DG_hu_h - D^2 u\|_{L^2(\mathcal{T}_h)}.
\end{equation*}

\begin{table}[htb!]
\centering
\caption{Numerical error of numerical test 1 on the uniform mesh}
\label{tab:test1_regular}
\resizebox{0.99\textwidth}{!}{
\begin{tabular}{|c|c|c|c|c|c|c|c|}
\hline 
 $\epsilon$ & $e$ & order& $De$ & order& $D^2e$ & order\\ \hline\hline 
  5.00e-01 &1.37e-01&--&6.36e-01&--&3.81e+00&--\\ \hline 
  2.50e-01 &1.06e-01&0.37&5.03e-01&0.34&3.37e+00&0.18\\ \hline 
  1.25e-01 &8.72e-02&0.28&4.22e-01&0.25&3.10e+00&0.12\\\hline 
  5.00e-02 &5.33e-02&0.54&2.75e-01&0.47&2.57e+00&0.20\\ \hline 
  2.50e-02 &3.19e-02&0.74&1.81e-01&0.60&2.19e+00&0.23\\ \hline 
  1.25e-02 &1.79e-02&0.84&1.15e-01&0.65&1.84e+00&0.25\\ \hline 
  5.00e-03 &7.84e-03&0.90&6.18e-02&0.68&1.45e+00&0.26\\ \hline 
  2.50e-03 &4.10e-03&0.93&3.84e-02&0.69&1.21e+00&0.27\\ \hline 
\end{tabular}}
\end{table}

\begin{table}[htb!]
\centering
\footnotesize
\caption{Numerical error of numerical test 1 on the unstructured mesh}
\label{tab:test1_del}
\resizebox{0.99\textwidth}{!}{
\begin{tabular}{|c|c|c|c|c|c|c|c|}
\hline 
 $\epsilon$ & $e$ & order& $De$ & order& $D^2e$ & order\\ \hline\hline 
 5.00e-01 &1.35e-01&--&6.24e-01&--&3.58e+00&--\\ \hline
 2.50e-01 &1.06e-01&0.35&4.97e-01&0.33&3.10e+00&0.21\\ \hline
 1.25e-01 &8.88e-02&0.26&4.19e-01&0.25&2.79e+00&0.15\\ \hline
 5.00e-02 &5.52e-02&0.52&2.71e-01&0.48&2.20e+00&0.26\\ \hline
 2.50e-02 &3.31e-02&0.74&1.73e-01&0.65&1.79e+00&0.30\\ \hline
 1.25e-02 &1.82e-02&0.86&1.03e-01&0.75&1.47e+00&0.29\\ \hline
 5.00e-03 &7.49e-03&0.97&4.85e-02&0.82&1.19e+00&0.23\\ \hline
 2.50e-03 &3.58e-03&1.07&3.05e-02&0.67&1.08e+00&0.14\\ \hline
\end{tabular}}
\end{table}

\subsection{Numerical Test 1} 

In this test, we choose the data $f(x,y) = 24x^2$ and $g(x,y) = x^4 + y^2$ to fit the unique convex classical solution $u(x,y) = x^4 + y^2$ of problem \eqref{equ: MA equation}. We choose the linear interpolation of $x^{1.2} + y^{1.2}$ as the initial guess $u_0$ for Newton’s method.

\begin{figure}[!h]
\captionsetup[subfigure]{labelformat=empty}
   \centering
   \subcaptionbox{$L^2$ error \label{fig:l2err_regular}}
  {\includegraphics[width=0.49\textwidth]{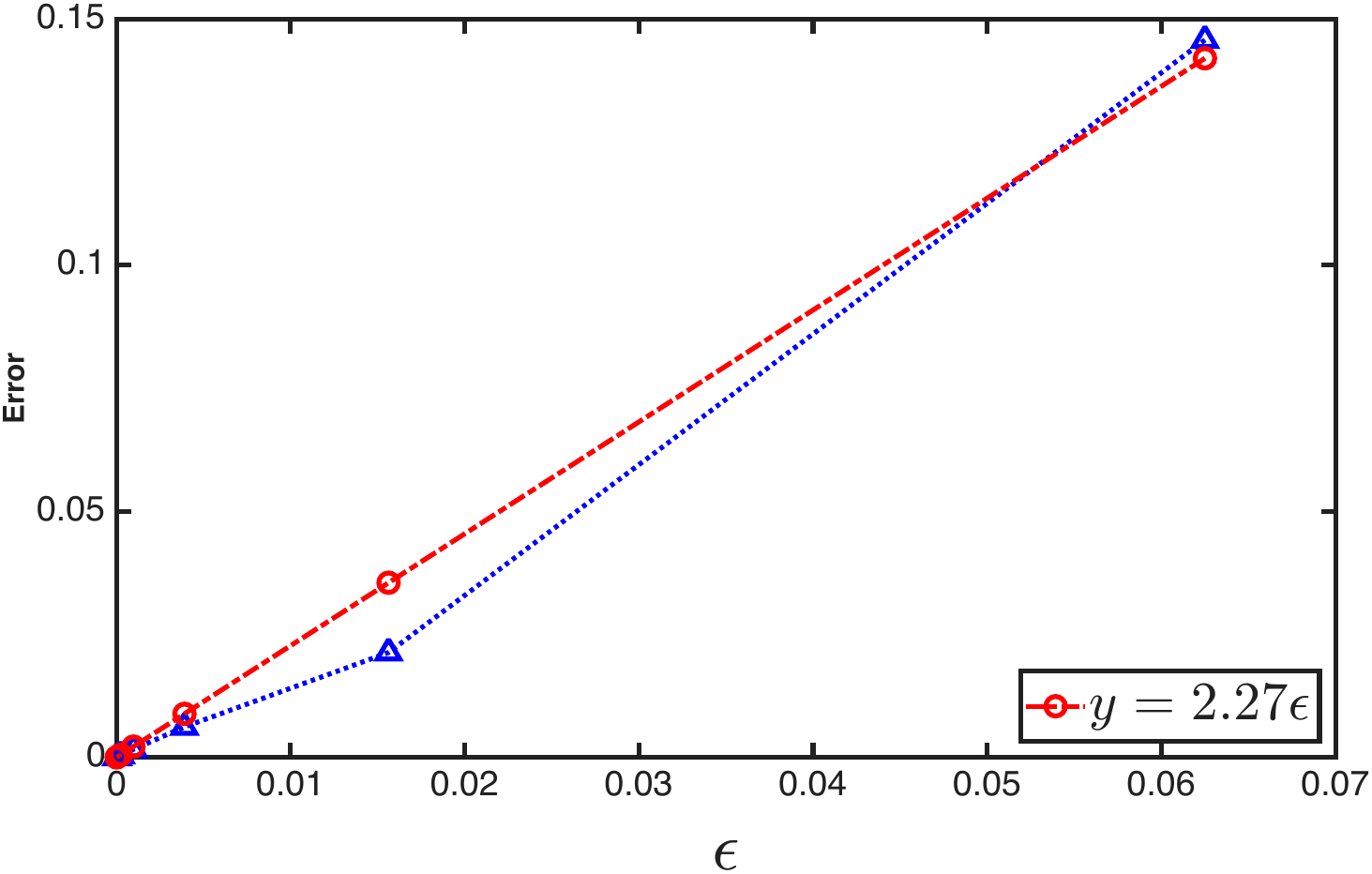}}
  \subcaptionbox{ $L^2$ error\label{fig:l2err_del}}
   {\includegraphics[width=0.49\textwidth]{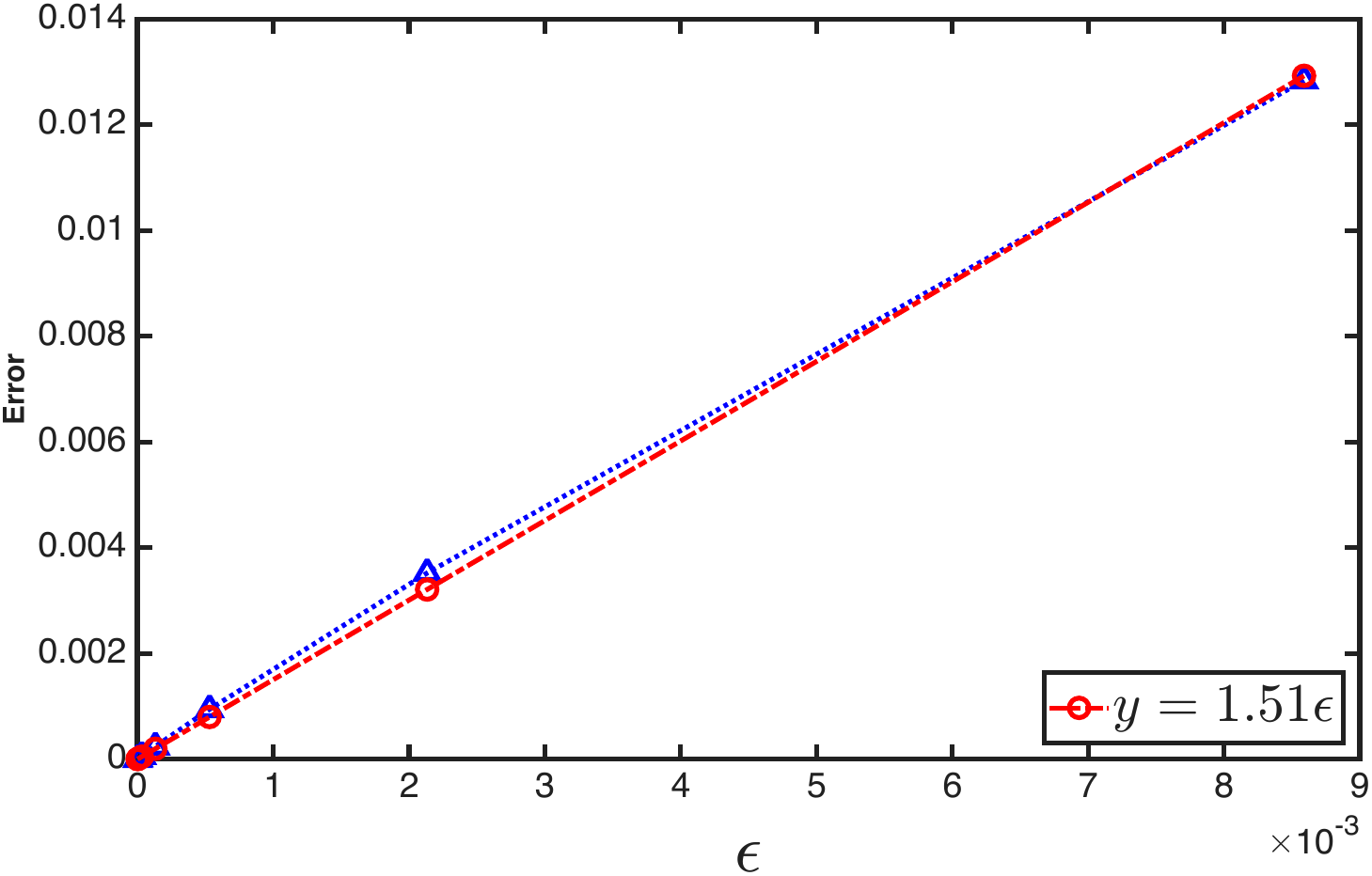}}
    \subcaptionbox{$H^1$ error \label{fig:h1err_regular}}
  {\includegraphics[width=0.49\textwidth]{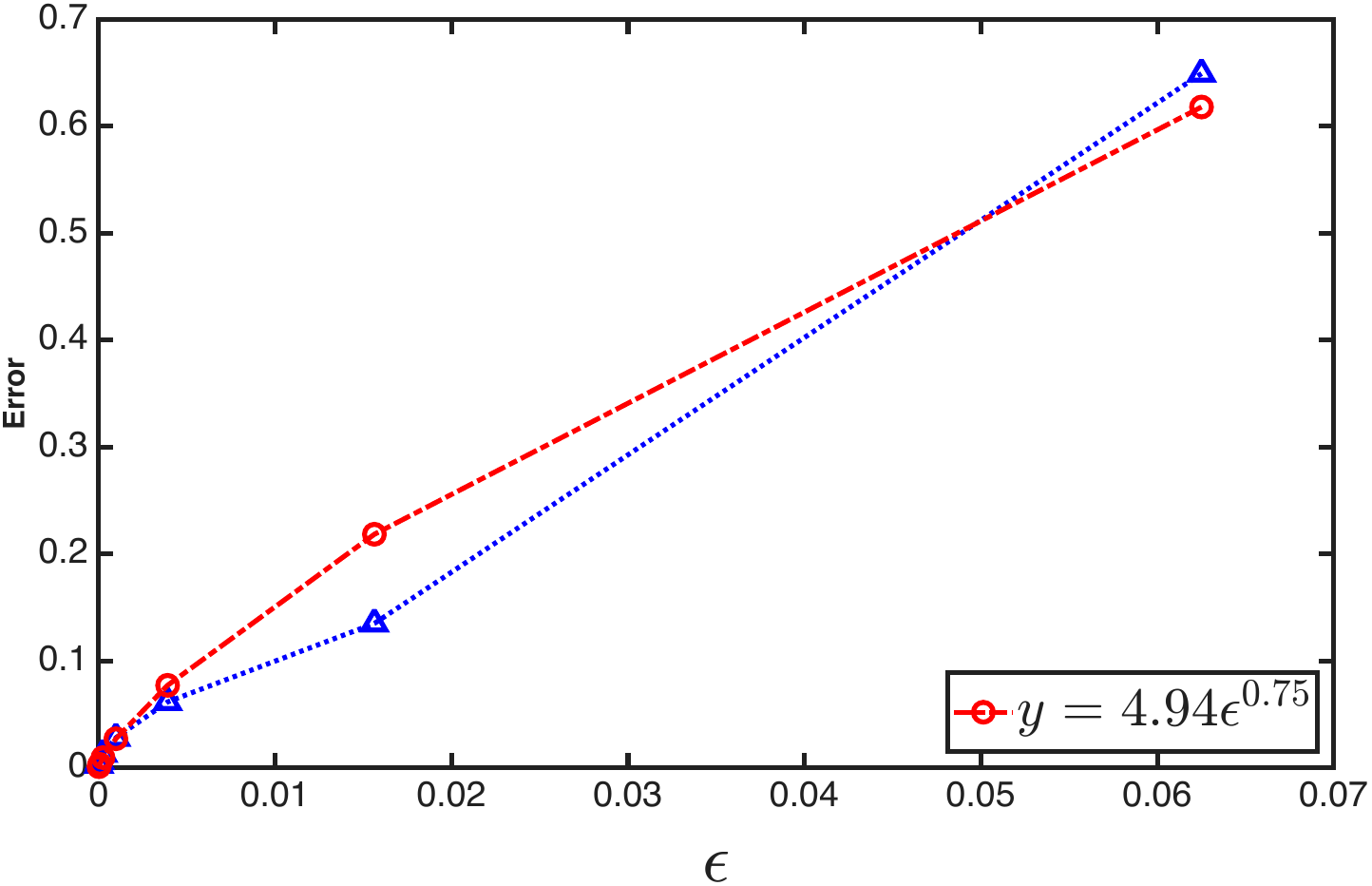}}
  \subcaptionbox{ $H^1$ error\label{fig:h1err_del}}
   {\includegraphics[width=0.49\textwidth]{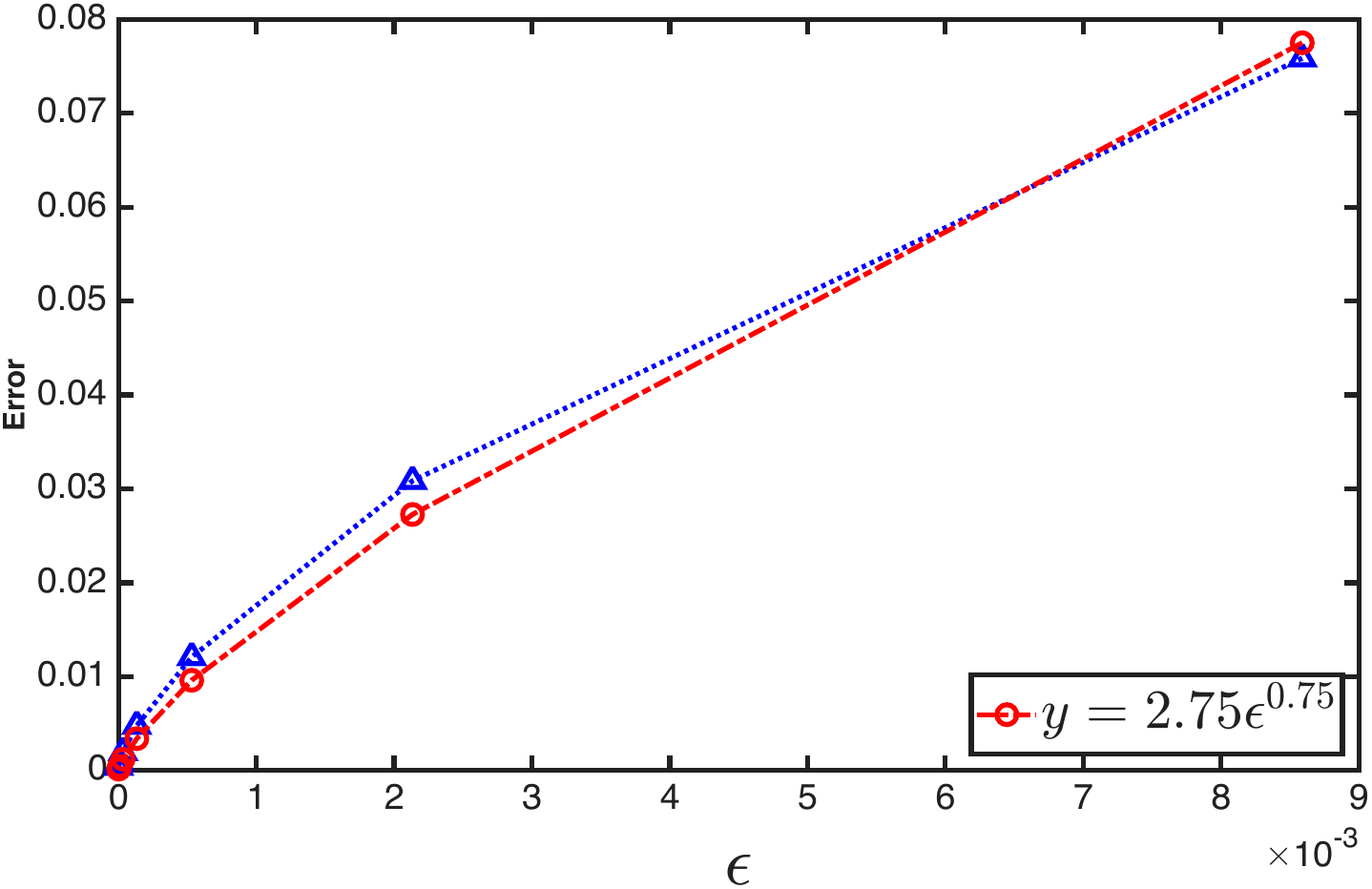}}
      \subcaptionbox{Discrete $H^2$ error \label{fig:h2err_regular}}
  {\includegraphics[width=0.49\textwidth]{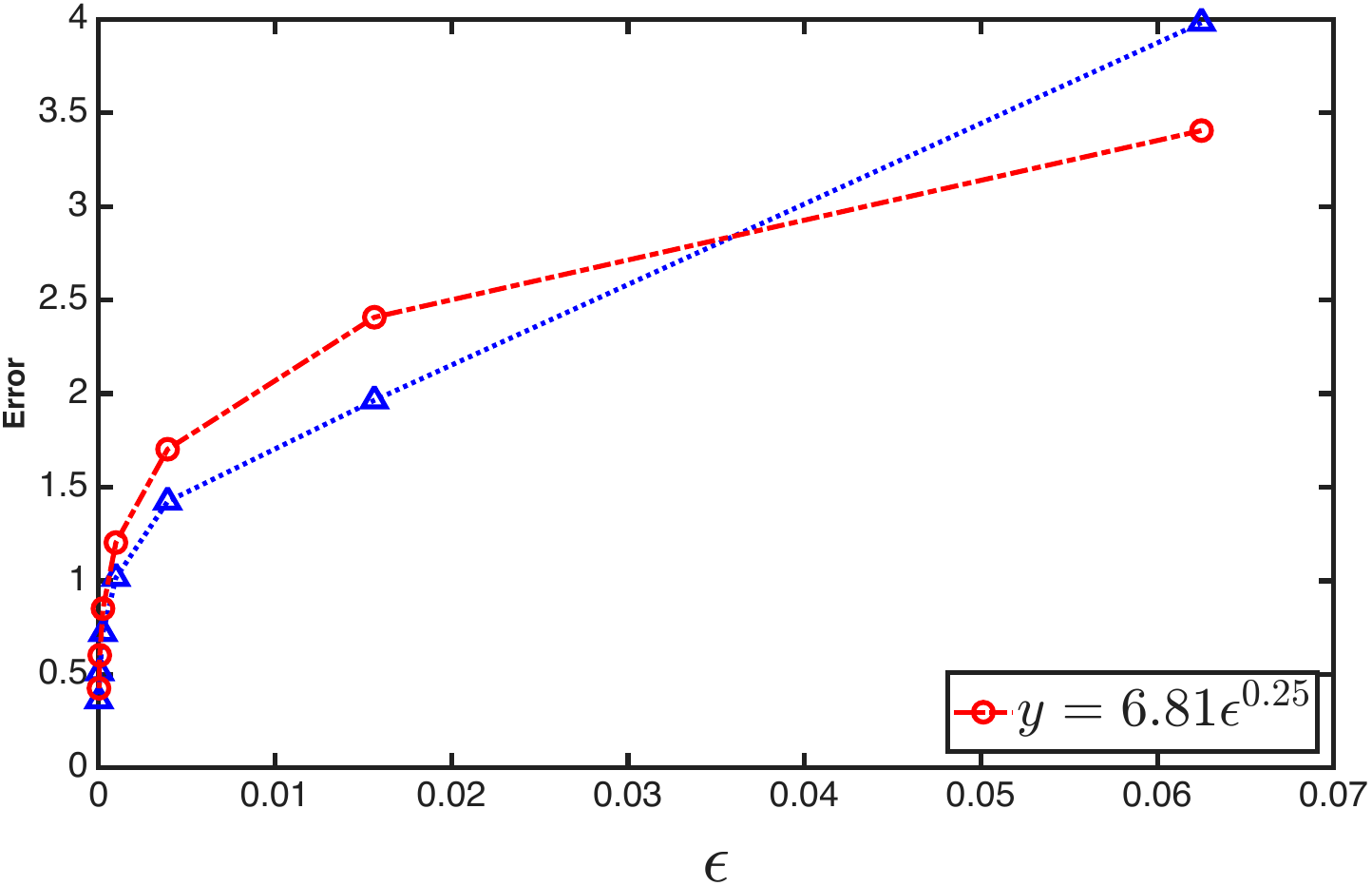}}
  \subcaptionbox{ Discrete $H^2$ error\label{fig:h2err_del}}
   {\includegraphics[width=0.49\textwidth]{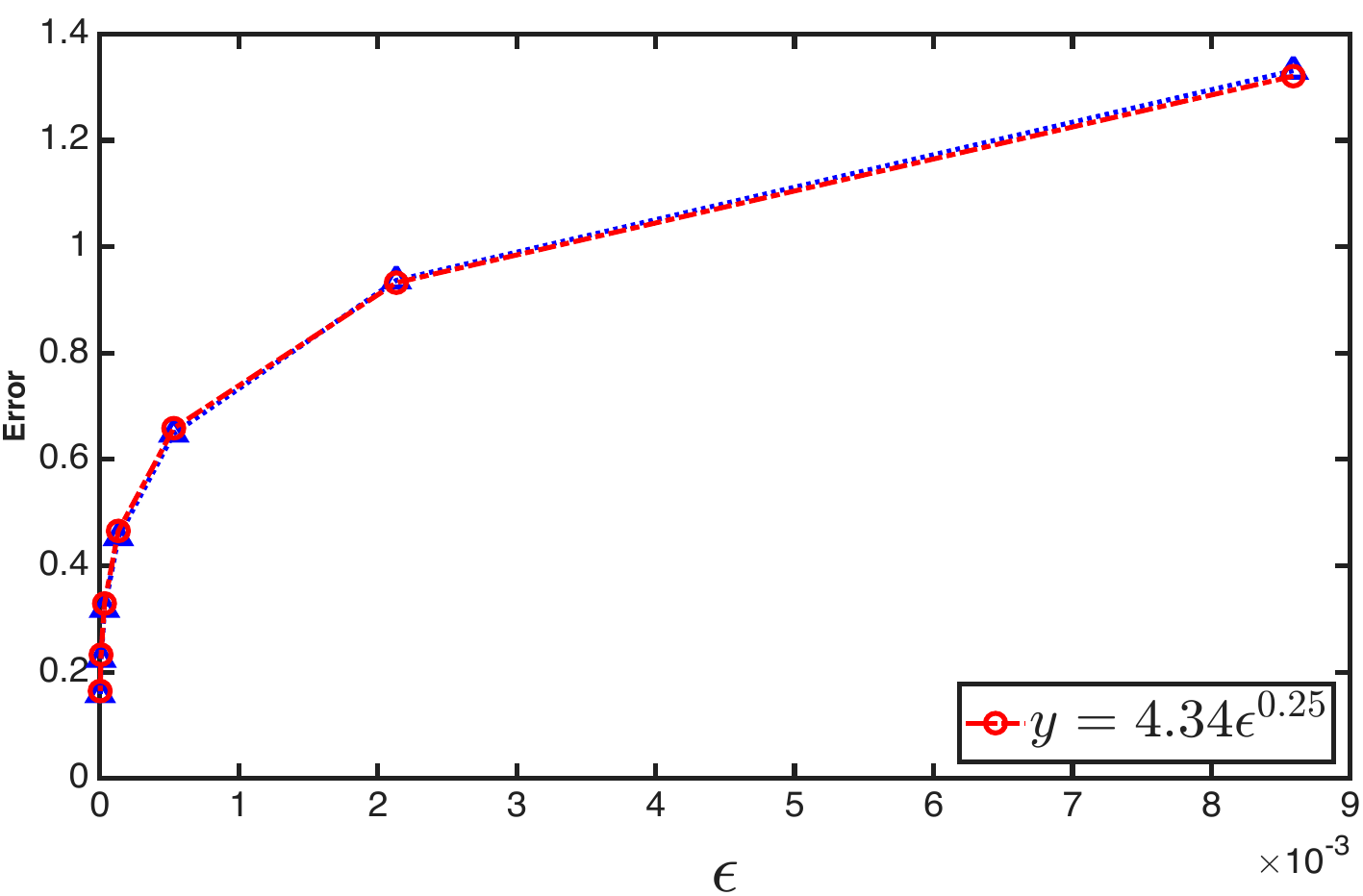}}
   \caption{Numerical errors with respect to $\epsilon$ and the fitted curves $y = \beta \eps^{\alpha}$.  Left: uniform meshes; right: unstructured meshes}\label{fig:errrel}
\end{figure}

We first focus on the numerical investigation of the behavior of the error $u - u_h^{\epsilon}$ in various norms, with a fixed small mesh size $h = 0.01$ and varying $\epsilon$. Given that $h$ is relatively small, $u_h^{\epsilon}$ is expected to be close to the regularized solution $u^{\epsilon}$. The numerical results are reported in Table~\ref{tab:test1_regular}. From the table, we clearly observe that the error $u - u_h^{\epsilon}$ decays at the rate of $\mathcal{O}(\epsilon)$ in the $L^2$-norm, $\mathcal{O}(\epsilon^{3/4})$ in the $H^1$-norm, and $\mathcal{O}(\epsilon^{1/4})$ in the discrete $H^2$-norm. This also suggests that the error $u - u^{\epsilon}$ decays at a similar rate. We also perform the test on an unstructured mesh, and the corresponding numerical errors are presented in Table~\ref{tab:test1_del}. Similar decay rates with respect to $\epsilon$ can be observed.

Then, we turn to explore a relation between $\epsilon$ and $h$ in order to determine an “optimal” choice of $h$ in terms of $\epsilon$ such that the global error $\left|u - u_h^\epsilon\right|$ exhibits the same convergence rate as $\left|u - u^\epsilon\right|$. We try several different relationships. The numerical results corresponding to the choice $h = \sqrt{\epsilon}$ are plotted in Figure~\ref{fig:errrel}. Based on the numerical data, we fit the constant $\beta$ in the expression $y = \beta \epsilon^{\alpha}$, where the exponent $\alpha$ varies depending on the norm. Specifically, we take $\alpha = 0.25$ for the discrete $H^2$ norm, $\alpha = 0.75$ for the $H^1$ norm, and $\alpha = 1$ for the $L^2$ norm, reflecting the distinct convergence behavior observed in each case. From the figure, we observe a close match between the numerical errors and the fitted curves, which provides evidence that $h = \sqrt{\epsilon}$ is indeed an “optimal” choice.

\subsection{Numerical Test 2}

The purpose of this test is to calculate the rate of convergence of $|u^{\epsilon} - u_h^{\epsilon}|$ for fixed $\epsilon = 0.001$ in various norms, where $u^{\epsilon}$ is the solution of \eqref{equ: model equations}. We select $f^{\epsilon}(x,y) = x^2 y^2 - 4\epsilon$, $g^{\epsilon}(x,y) = \frac{1}{12}(x^4 + y^4)$, and $\psi^{\epsilon}(x,y) = x^2 + y^2$ so that the solution is given by $u^{\epsilon} = \frac{1}{12}(x^4 + y^4)$.

The numerical results on uniform meshes are summarized in Table~\ref{tab:test2_regular}. It can be seen from the data that the discrete $H^2$ error converges at the rate of $\mathcal{O}(h^{0.5})$, the $H^1$ error at the rate of $\mathcal{O}(h)$, and the $L^2$ error at the rate of $\mathcal{O}(h^2)$. Compared to the existing recovery-based finite element method for the Monge–Ampère equation \cite{chenfengzhang2021}, our method achieves improved convergence rates in terms of both $L^2$ and $H^1$ errors. The numerical results on unstructured meshes are reported in Table~\ref{tab:test2_del}, where similar convergence rates can be observed for the $H^1$ and discrete $H^2$ errors.

\begin{table}[htb!]
\centering
\footnotesize
\caption{Numerical error of numerical test 2 on uniform meshes}
\label{tab:test2_regular}
\resizebox{0.99\textwidth}{!}{
\begin{tabular}{|c|c|c|c|c|c|c|c|}
\hline 
 $1/h$ & $e$ & order& $De$ & order& $D^2e$ & order\\ \hline\hline 
  1.60e+01 &1.63e-04&--&6.85e-03&--&2.45e-01&--\\  \hline
  3.20e+01 &4.64e-05&1.82&3.37e-03&1.03&1.76e-01&0.48\\ \hline
  6.40e+01 &1.29e-05&1.84&1.67e-03&1.01&1.26e-01&0.49\\ \hline
  1.28e+02 &3.47e-06&1.90&8.32e-04&1.01&8.94e-02&0.49\\ \hline
  2.56e+02 &9.02e-07&1.94&4.14e-04&1.01&6.34e-02&0.50\\ \hline
\end{tabular}}
\end{table}

\begin{table}[htb!]
\centering
\footnotesize
\caption{Numerical error of numerical test 2 on unstructured meshes}
\label{tab:test2_del}
\resizebox{0.99\textwidth}{!}{
\begin{tabular}{|c|c|c|c|c|c|c|c|}
\hline 
 $1/h$ & $e$ & order& $De$ & order& $D^2e$ & order\\ \hline\hline 
 1.08e+01 &3.63e-04&--&7.53e-03&--&2.70e-01&--\\  \hline
 2.16e+01 &7.60e-05&2.24&3.39e-03&1.15&1.93e-01&0.48\\  \hline
 4.34e+01 &2.29e-05&1.73&1.65e-03&1.04&1.38e-01&0.48\\  \hline
 8.68e+01 &7.62e-06&1.58&8.17e-04&1.01&9.85e-02&0.49\\  \hline
 1.74e+02 &2.95e-06&1.37&4.06e-04&1.01&6.99e-02&0.49\\  \hline
\end{tabular}}
\end{table}

\subsection{Numerical Test 3}

\begin{figure}[!h]
   \centering
  \includegraphics[width=0.65\textwidth]{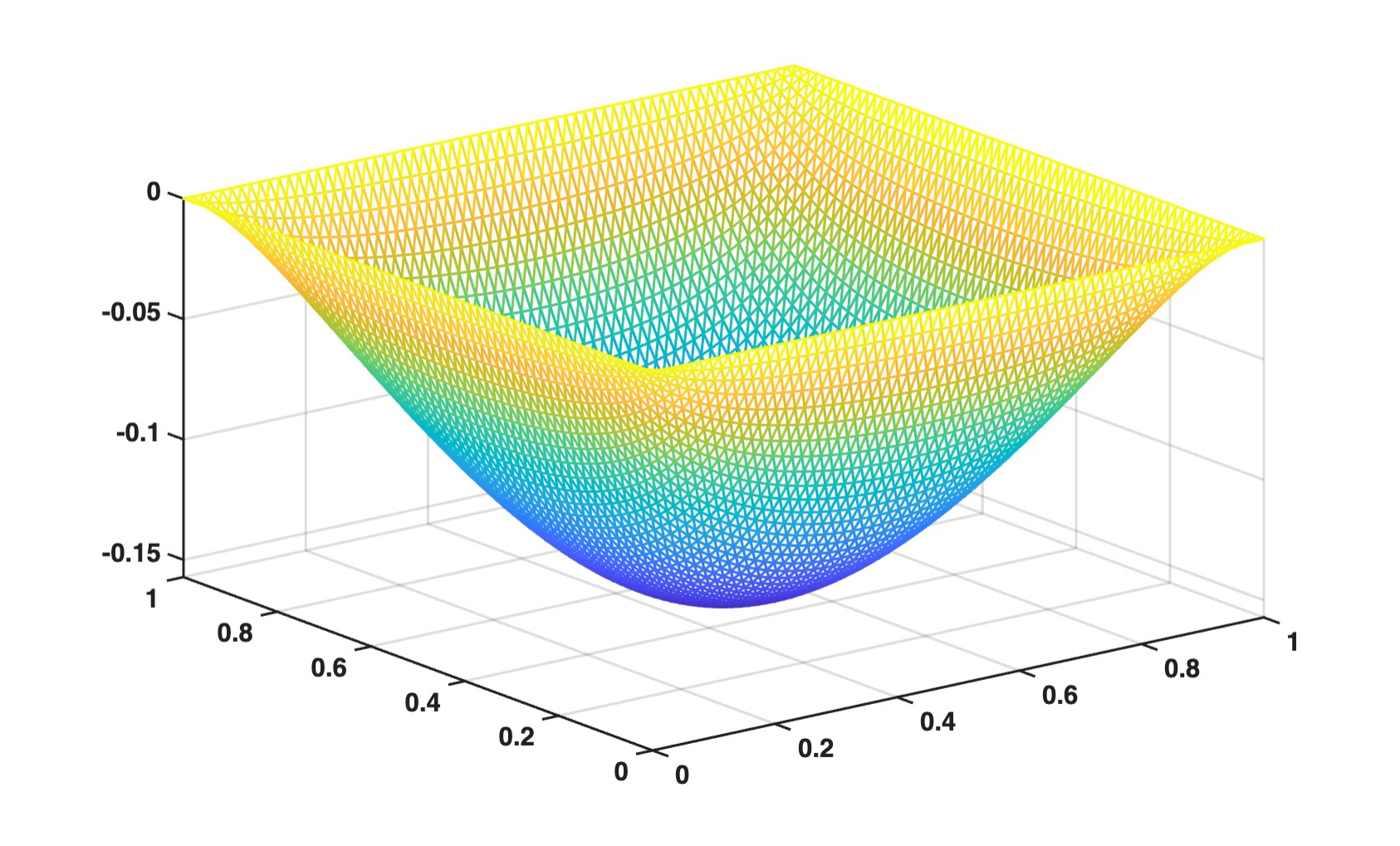}
   \caption{Numerical solution of Test 3 computed by the proposed linear element method with $\epsilon=0.001$ and $h = 1/64$.}
   \label{fig:vissol}
\end{figure}

In this test, we choose $f = 1$ and $g = 0$ as the data for problem \eqref{equ: MA equation}. We note that problem \eqref{equ: MA equation} admits a unique convex viscosity solution but does not possess a classical solution. We impose the boundary condition $\psi^{\epsilon} \equiv \epsilon$. The numerical solution is displayed in Figure~\ref{fig:vissol}. From the figure, it is clear that the computed solution successfully captures the convexity of the nonsmooth solution.

\subsection{Numerical Test 4} 

In this experiment, we consider the computational domain $\Omega = (-1,1)^2$. The exact solution is $u(x,y)=|x| \in H^1(\Omega)$, which satisfies the Monge–Amp\`ere equation in the viscosity sense with $f=0$. Because of the extremely low regularity of $u$, no convergence rate can be expected. In the numerical tests, we set $\epsilon = 0.001$, and the resulting errors are summarized in Table~\ref{tab:singular}. The numerical results shows the error is decreasing.

\begin{table}[ht]
\centering
\caption{Numerical results for the Monge--Amp\`ere equation with solution $u=|x|$.}
\label{tab:singular}
\begin{tabular}{c c c c c}
\hline
$N$ & $32$ & $64$ & $128$ & $256$ \\
\hline
$L^2$ error & $1.16\times10^{-1}$ & $9.00\times10^{-2}$ & $8.11\times10^{-2}$ & $7.89\times10^{-2}$ \\
Order & -- & $0.37$ & $0.15$ & $0.04$ \\
\hline
\end{tabular}
\end{table}

\section{Conclusion} \label{sec:con}

In this paper, we propose a recovery-based linear finite element method for solving the regularized Monge--Amp\`ere equation, inspired by the discrete MT estimate of the gradient recovery operator. We establish convergence in both the discrete $H^2$ and $H^1$ norms. Numerical results show that the proposed method outperforms the existing recovery-based finite element method for the Monge-Amp\`ere equation in terms of accuracy and convergence rates.

\section*{Acknowledgments}

This work was supported in part by the Dyason Fellowship and the Faculty Science Researcher Development Grant of the University of Melbourne, and the NSFC grant 12131005. 

\appendix

\section{Properties of the recovery operator}\label{sec: appendix}

In this appendix, we list some properties of $G_h$. Proofs of these properties can be found in \cite{chu2025recovery}.

\begin{lemma}\label{lem: apo of function in V_h}
	For any $v_h \in V_h$, we have 
	\begin{equation*}
	\|\nabla v_h - G_hv_h\|_{L^2(\Omega)} \le C h\|DG_hv_h\|_{L^2(\mathcal{T}_h)}.
	\end{equation*}
	Consequently, we obtain the following bound:
	\begin{equation*}
	\|G_hv_h\|_{L^2(\Omega)} \le C\|\nabla v_h\|_{L^2(\Omega)} \le C \|G_hv_h\|_{L^2(\Omega)}.
	\end{equation*}
\end{lemma}

\begin{lemma}\label{lem: dis poin 2}
	For any $v_h \in V_h^0$, we have 
	\begin{equation*}
	\|\nabla v_h\|_{L^2(\Omega)} \le C  \|DG_hv_h\|_{L^2(\mathcal{T}_h)}.
	\end{equation*}
\end{lemma}

\begin{theorem}\label{thm: AP on a_mesh}
	Assuming that $v \in W^{2,p}(\Omega)$ with $1 \le p \le \infty$. Then the following local approximation property holds:
	\begin{equation*}
	\|\nabla v - G_hv_I \|_{L^p(T)} \le C_p  h_T |v|_{W^{2,p}(\mathcal{K}_T)}, \qquad \forall T \in \mathcal{T}_h.
	\end{equation*}
	This yields a global approximation property:
	\begin{equation*}
	\|\nabla v - G_hv_I \|_{L^p(\Omega)} \le C_p h |v|_{W^{2,p}(\Omega)}.
	\end{equation*}
\end{theorem}

\begin{theorem}\label{thm: Approximation properties on general mesh}
	Let $\mathcal{T}_h$ be a mesh that satisfies the $(\nu, \mu)$ condition (see \cite{chu2025recovery}). Let $\rho=\min\{1,\nu,\mu/2\}$, and let $v\in H^3(\Omega)$ and $h\leq e^{-2}$. Then, for $\beta \in [0,2]$, we have the following global superconvergence result:
	\begin{equation*}
	\left(\sum_{T \in \mathcal{T}_h} h_T^{-\beta}\|\nabla v - G_hI_hv \|_{L^2(T)}^2 \right)^{1/2} \le C h^{(1 + \rho) - \beta/2} |\log h|^{1/2} \|v\|_{H^3(\Omega)}.
	\end{equation*}
\end{theorem}

Next, we will prove a superconvergence result of the operator $G_h$.

\begin{theorem}\label{thm: superconvergence property of G_h}
    Suppose ${\bm v} \in \Big(H^1(\Omega)\Big)^2  $ and $v_h \in V_h$. Then 
    \begin{equation*}
        \Big|\int_{\Omega} {\bm v} \cdot (\nabla v_h - G_hv_h ) \ \text{d}x \Big| \le C h^2 |{\bm v}|_{H^1(\Omega)} \|DG_hv_h\|_{L^2(\mathcal{T}_h)}.
    \end{equation*}
\end{theorem}
\begin{proof}
    Let $\{\phi_m \}_{m \in \mathcal{M}_h}$ be the nodal basis of $V_h^{CR}$, and let $\T_m$ be the union of the elements sharing the midpoint $m$. Then for any ${\bm w}_h \in V_h^{CR}\times V_h^{CR}$, we have 
    \begin{align*}
&\int_{\Omega} \sum_{m \in \mathcal{M}_h}{\bm w}_h(m)\cdot G_hv_h(m) \phi_m \ \text{d}x \\
= & \sum_{m \in \mathcal{M}_h} {\bm w}_h(m)\cdot G_hv_h(m) \frac{|\mathcal{T}_m|}{3} 
= \sum_{m \in \mathcal{M}_h} \frac{{\bm w}_h(m)}{3}\cdot \int_{\mathcal{T}_m} \nabla v_h \ \text{d}x\\
=& \sum_{T \in \mathcal{T}_h} \int_{T}\Big(  \sum_{m \in \mathcal{M}_h \cap T} \frac{{\bm w}_h(m)}{3} \Big) \cdot \nabla v_h \ \text{d}x 
 = \sum_{T \in \mathcal{T}_h} \int_{T} \sum_{m \in \mathcal{M}_h \cap T} {\bm w}_h(m)\phi_m\cdot \nabla v_h \ \text{d}x \\
= &\sum_{T \in \mathcal{T}_h} \int_{T} {\bm w}_h\cdot \nabla v_h = \int_{\Omega} {\bm w}_h \cdot \nabla v_h \ \text{d}x.
\end{align*}
where we have used the fact that $\nabla v_h$ is constant on $T$. Let ${\bm w}_h = ({\bm v})_{I^{CR}}$ which denotes the interpolation on the space $V_h^{CR}$, then we have 
\begin{align*}
    \| {\bm v} - {\bm w}_h \|_{L^2(\Omega)} & \le C h |{\bm v}|_{H^1(\Omega)}, \\
    \| D {\bm w}_h \|_{L^2(\mathcal{T}_h)} & \le C |{\bm v}|_{H^1(\Omega)}. 
\end{align*}

We have the decomposition 
\begin{align*}
    &\quad\int_{\Omega} {\bm v} \cdot (\nabla v_h - G_hv_h ) \ \text{d}x\\
    &=  \underbrace{\int_{\Omega}( {\bm v} - {\bm w}_h) \cdot (\nabla v_h - G_hv_h ) \ \text{d}x}_{=:I_1} + \underbrace{\int_{\Omega} {\bm w}_h \cdot (\nabla v_h - G_hv_h ) \ \text{d}x}_{=:I_2}.
\end{align*}
For $I_1$, we have 
\begin{equation*}
    |I_1|\le \|{\bm v} - {\bm w}_h\|_{L^2(\Omega)}\|\nabla v_h - G_hv_h\|_{L^2(\Omega)} \le C h^2 |{\bm v}|_{H^1(\Omega)} \|DG_hv_h\|_{L^2(\mathcal{T}_h)},
\end{equation*}
and for $I_2$, we have 
\begin{align*}
    I_2 & = \int_{\Omega} {\bm w}_h \cdot (\nabla v_h - G_hv_h ) \ \text{d}x \\
    & = \int_{\Omega} {\bm w}_h \cdot \nabla v_h  \ \text{d}x - \int_{\Omega} {\bm w}_h \cdot  G_hv_h  \ \text{d}x \\
    & = \int_{\Omega} \sum_{m \in \mathcal{M}_h}{\bm w}_h(m)\cdot G_hv_h(m) \phi_m \ \text{d}x - \int_{\Omega} {\bm w}_h \cdot  G_hv_h  \ \text{d}x \\
    & = \sum_{T \in \mathcal{T}_h} \int_T \Big( \big(  {\bm w}_h \cdot  G_hv_h \big)_{I^{CR}} - {\bm w}_h \cdot  G_hv_h \Big) \ \text{d}x,
\end{align*}
then we have 
\begin{align*}
    |I_2| & \le \sum_{T \in \mathcal{T}_h} \Big\|  \big(  {\bm w}_h \cdot  G_hv_h \big)_{I^{CR}} - {\bm w}_h \cdot  G_hv_h  \Big\|_{L^1(T)}  \le C h^2 \sum_{T \in \mathcal{T}_h} \|D^2({\bm w}_h \cdot  G_hv_h ) \|_{L^1(T)} \\
    & \le Ch^2 \sum_{T \in \mathcal{T}_h} \|D{\bm w}_h\|_{L^2(T)}\|DG_hv_h\|_{L^2(T)}  \le Ch^2 \| D {\bm w}_h \|_{L^2(\mathcal{T}_h)} \| D G_h v_h \|_{L^2(\mathcal{T}_h)} \\
    &\le C h^2 |{\bm v}|_{H^1(\Omega)} \|DG_hv_h\|_{L^2(\mathcal{T}_h)}.
\end{align*}
Then we can complete the proof by the above estimates.
\end{proof}

\begin{lemma}\label{lem: discrete infty}
    There holds
    \begin{align}
        \|v_h\|_{L^{\infty}(\Omega)} \le C\|\Delta_h v_h\|_{L^2(\Omega)} \qquad \forall\, v_h \in V_h^0,
    \end{align}
    where $\Delta_h: V_h^0 \to V_h^0$ is the discrete Laplacian operator defined as 
    \begin{align}
        (\Delta_h u_h, v_h)_{L^2(\Omega)} = -(\nabla u_h, \nabla v_h)_{L^2(\Omega)} \qquad \forall\, v_h \in V_h^0.
    \end{align}
\end{lemma}
\begin{proof}
    By the discrete Meyers estimate in \cite{badr2014meyers}, there exists $\varepsilon > 0$ such that for every $p \in (2,2+\varepsilon)$, 
    \begin{align*}
    	\|v_h\|_{W^{1,p}(\Omega)} \le C_p \|\Delta_h v_h\|_{W^{-1,p}(\Omega)}.
    \end{align*}
    Fix such a $p$ and let $p' = p/(p-1)$, there holds $W_0^{1,p'}(\Omega) \hookrightarrow L^2(\Omega)$. For any $\varphi\in W_0^{1,p'}(\Omega)$, we have 
    \begin{align*}
    	|(\Delta_h v_h, \varphi)_{L^2(\Omega)}| \le \|\Delta_h v_h\|_{L^2(\Omega)}\|\varphi\|_{L^2(\Omega)} \le C \|\Delta_h v_h\|_{L^2(\Omega)}\|\varphi\|_{W_0^{1,p'}(\Omega)}.
    \end{align*}
    This implies 
    \begin{align*}
    	\|\Delta_h v_h \|_{W^{-1,p}(\Omega)} \le C \|\Delta_h v_h\|_{L^2(\Omega)}.
    \end{align*}
    Finally, by the embedding $W^{1,p}(\Omega) \hookrightarrow L^{\infty}(\Omega)$, we obtain 
    \begin{align*}
    	\|v_h\|_{L^{\infty}(\Omega)} \le C \|v_h\|_{W^{1,p}(\Omega)} \le C \|\Delta_h v_h \|_{W^{-1,p}(\Omega)} \le C \|\Delta_h v_h\|_{L^2(\Omega)},
    \end{align*}
    which completes the proof.
    \end{proof}
\begin{lemma}\label{lem: discrete Poincare inequality}
    We have the following estimate:
    \begin{align}
        \|v_h\|_{L^{\infty}(\Omega)} \le C \|DG_h v_h\|_{L^2(\T_h)} \qquad \forall\, v_h \in V_h^0.
    \end{align}
\end{lemma}
\begin{proof}
    From Lemma \ref{lem: discrete infty}, we only need to prove that 
    \begin{align}
        \|\Delta_h v_h\|_{L^2(\Omega)} \le C \|DG_h v_h\|_{L^2(\T_h)}.
    \end{align}
    Note that $\|\Delta_h v_h\|_{L^2(\Omega)} = \sup_{0 \neq u_h \in V_h} |(\Delta_h v_h, u_h)_{L^2(\Omega)}|/\|u_h\|_{L^2(\Omega)} $, then fix $0 \neq u_h \in V_h$, 
    \begin{align*}
        &(\Delta_h v_h, u_h)_{L^2(\Omega)} \\
        & = - (\nabla v_h, \nabla u_h)_{L^2(\Omega)} \\
        & = (G_hv_h - \nabla v_h,  \nabla u_h)_{L^2(\Omega)} - (G_hv_h, \nabla u_h)_{L^2(\Omega)} \\
        & = (G_hv_h - \nabla v_h,  \nabla u_h)_{L^2(\Omega)} + (\text{div}(G_hv_h), u_h)_{L^2(\Omega)}  - \sum_{e \in \mathcal{E}_h^I} ( \lj G_hv_h \rj\cdot \bm{n}_e, u_h)_{L^2(e)}.
    \end{align*}
    The above identity implies 
    \begin{align*}
        |(\Delta_h v_h, u_h)_{L^2(\Omega)}| & \le Ch \|DG_hv_h\|_{L^2(\T_h)} \|\nabla u_h\|_{L^2(\Omega)} +  \|DG_hv_h\|_{L^2(\T_h)} \|u_h\|_{L^2(\Omega)} \\
        & \quad + \left( \sum_{e \in \mathcal{E}_h^I} h_e^{-1} \| \lj G_hv_h \rj\|_{L^2(e)}^2 \right)^{\frac12}\left( \sum_{e \in \mathcal{E}_h^I} h_e \|  u_h \|_{L^2(e)}^2 \right)^{\frac12} \\
        & \le  C \|DG_hv_h\|_{L^2(\T_h)} \|u_h\|_{L^2(\Omega)}.
    \end{align*}
    Here we used the inverse inequality and trace inequality in the last line. This completes the proof.
\end{proof}

\bibliographystyle{siamplain}
\bibliography{references}

\end{document}